\documentclass[11pt,leqno,twoside]{amsart}
\usepackage{amssymb,verbatim,graphicx}

\renewcommand{\baselinestretch}{1.05} 

\def\figdir{figs}
\def\incgr#1#2{\includegraphics[width=#1]{\figdir/#2}}
\def\trisph#1{\vbox{
    \hbox{\incgr{45mm}{tri#1}}\vspace{5mm}
          \hbox{\hspace{5mm}\incgr{35mm}{sph#1}}\vspace{2mm} }}
\def\figsx{\hbox{\trisph{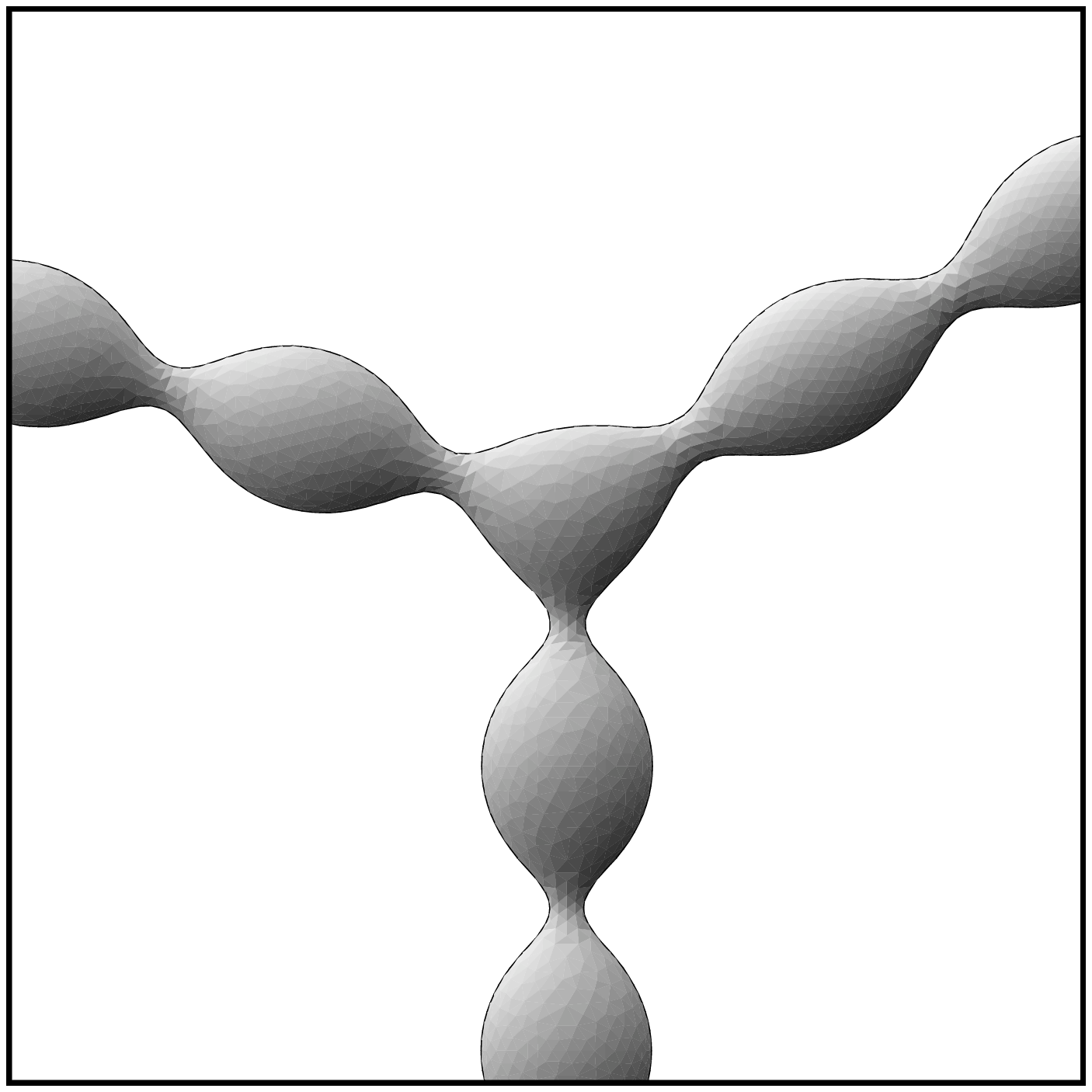}\hspace{5mm}\trisph{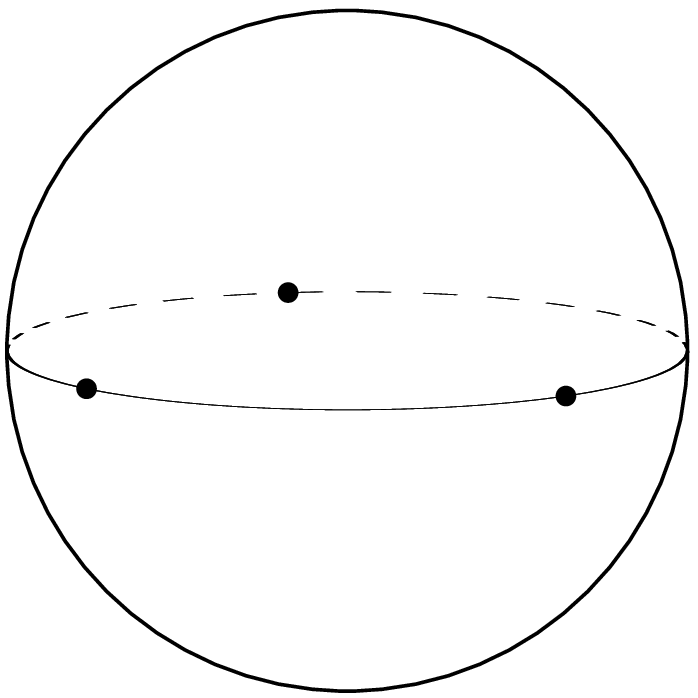}\hspace{5mm}\trisph{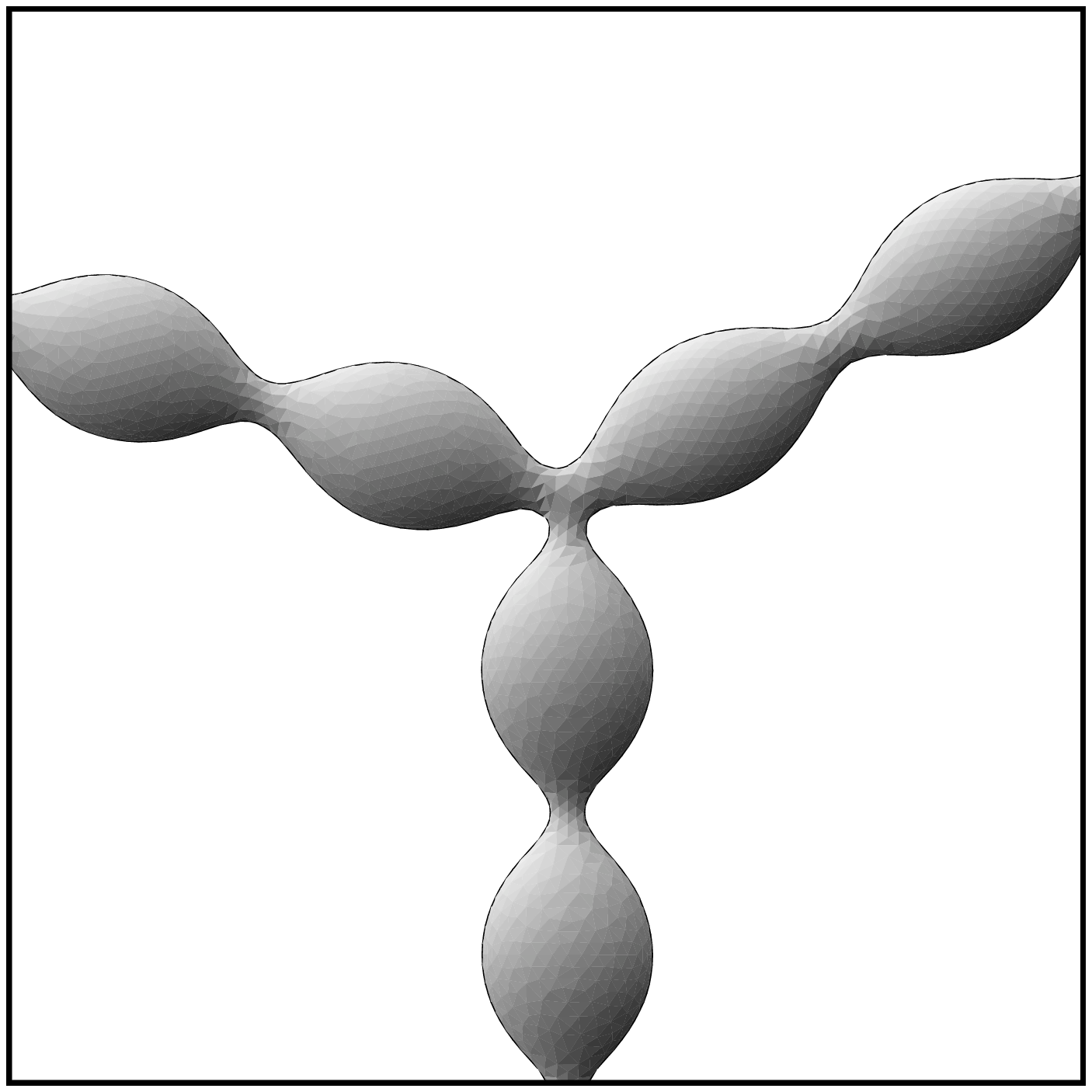}}}
\def\figsix#1#2{\figfig{triunds}{\figsx}{#1}{#2}}
\def\figtwo#1#2{\figfig{covsp}
    {\incgr{60mm}{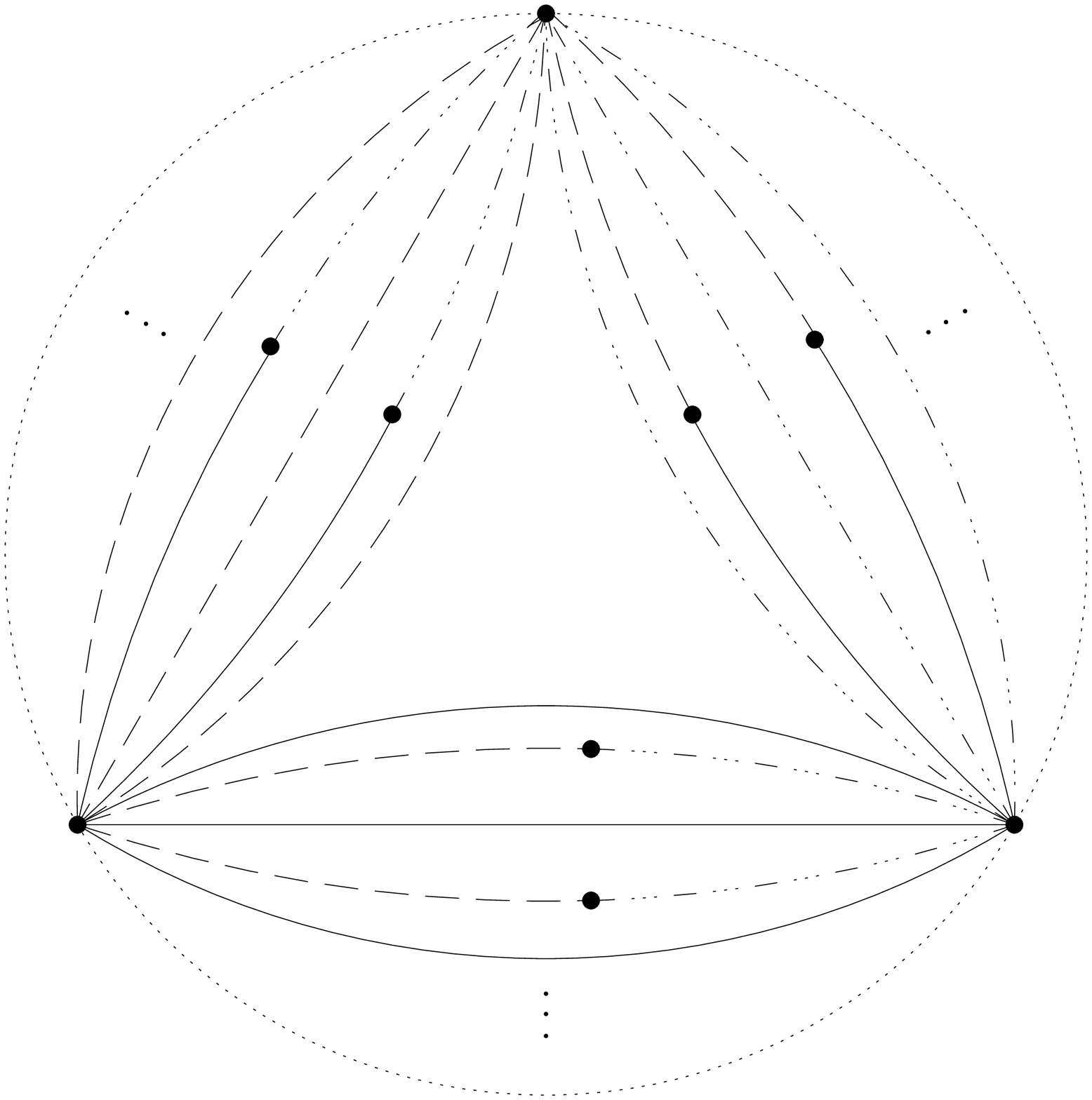}\hspace{15mm}\incgr{70mm}{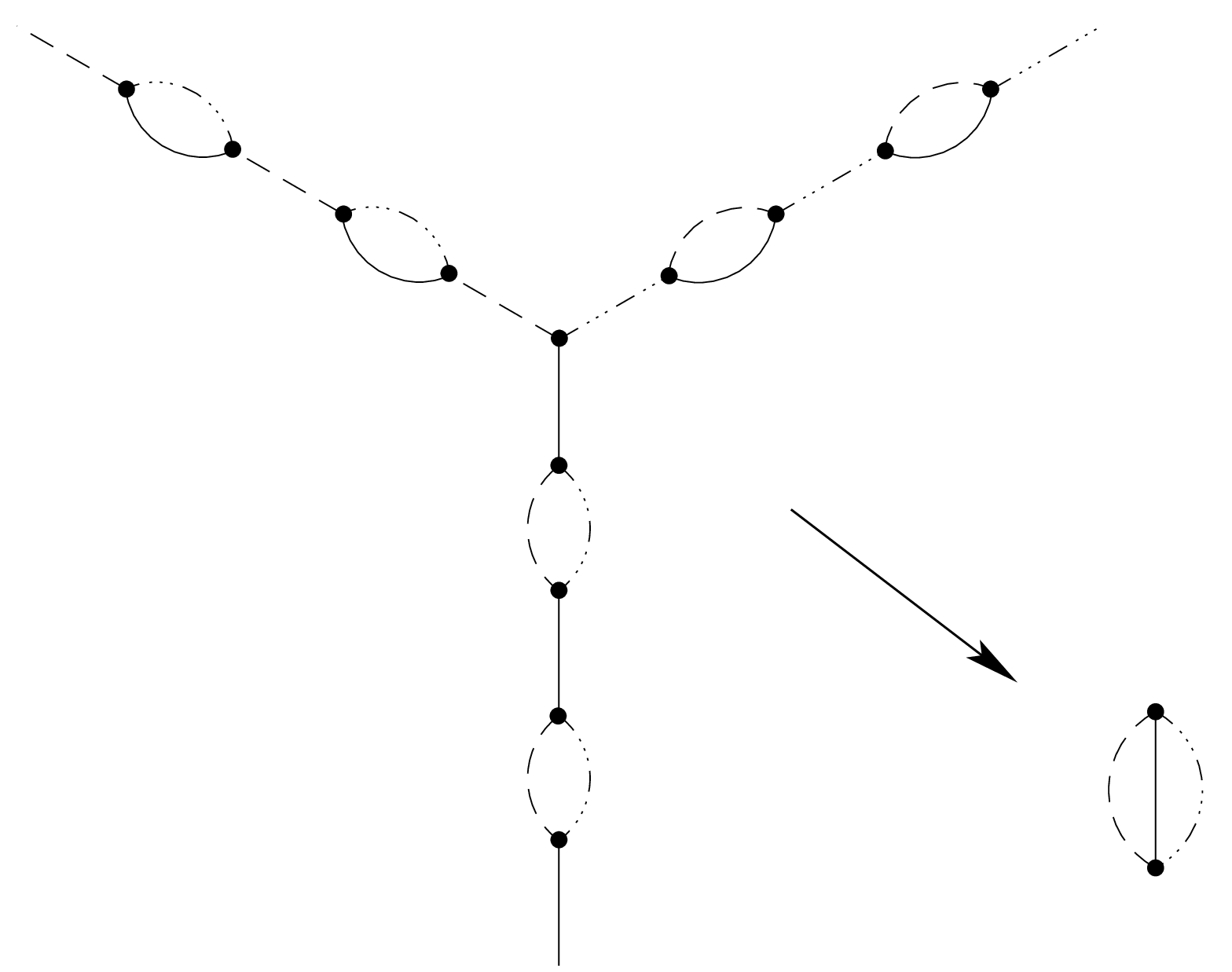}} {#1}{#2}}

\def\figfig#1#2#3#4{\begin{figure}
    \centerline{#2} \caption[#3]{{\small #4}}\label{jfig:#1}
    \end{figure}}

\def\figr#1{Figure~\ref{jfig:#1}}

\newtheorem{theorem}{Theorem}[section]
\newtheorem*{maintheorem}{Main Theorem}
\newtheorem{lemma}[theorem]{Lemma}
\newtheorem{proposition}[theorem]{Proposition}
\newtheorem{corollary}[theorem]{Corollary}
\newtheorem*{cornn}{Corollary}
\theoremstyle{definition}
  \newtheorem*{definition}{Definition}         
\theoremstyle{remark}
  \newtheorem*{remark}{Remark}
\def\eqn#1{\eqref{eq:#1}} 

\newcommand{\area}{\operatorname{area}}
\def\alex-emb/{Alexandrov-embedded}
\newcommand{\after}{\circ}
\renewcommand{\bar}{\overline}

\def\CMC/{{\Small{CMC}}}
\def\d{\partial}

\renewcommand{\epsilon}{\varepsilon}
\renewcommand{\H}{\mathbb{H}}
\let\dotlessi\i
\renewcommand{\i}{\mathbf{i}}
\newcommand{\id}{\operatorname{Id}}
\newcommand{\ident}{\equiv}
\newcommand{\isom}{\cong}
\renewcommand{\Im}{\operatorname{Im}}
\renewcommand{\j}{\mathbf{j}}
\renewcommand{\k}{\mathbf{k}}
\def\mod#1{\bmod{\,#1}}
\newcommand{\M}{\mathcal{M}}
\renewcommand{\phi}{\varphi}
\newcommand{\Pik}{\Pi_{\k}}               
\newcommand{\R}{\mathbb{R}}
\renewcommand{\Re}{\operatorname{Re}}
\let\paragraph=\S  
\renewcommand{\S}{\mathbb{S}}
\newcommand{\SO}{\mathsf{SO}}

\def\setm{\smallsetminus}
\newcommand{\T}{\mathcal{T}}
\def\threept/{three-point}
\def\takes{\colon}

\def\kHopf/{$\k$-Hopf}

\def\wtilde#1{{\mskip4mu\widetilde{\mskip-4mu #1\mskip-1mu}\mskip1mu}}
\def\plustilde#1{\wtilde{#1}{}^+}
\def\Mconj{{\wtilde M}}
\def\Sconj{{\wtilde S}}
\def\fconj{{\tilde f}}
\def\dfconj{{d\!\tilde f}}
\def\gammaconj{{\tilde \gamma}}

\def\nuconj{{\tilde \nu}}
\def\etaconj{{\tilde \eta}}

\begin{document}

\title[Triunduloids]
{Triunduloids:\\ Embedded constant mean curvature surfaces\\
with three ends and genus zero}
\date{January 2001; minor revisions July 2002 and April 2003}

\author[Grosse-Brauckmann]{Karsten Gro\ss e-Brauckmann}
\address{Technische Universit\"at Darmstadt, Fachbereich 
    Mathematik (AG 3), Schlossgartenstr. 7, 64289 Darmstadt, Germany}
    \email{kgb@math.uni-bonn.de}

\author[Kusner]{Robert B. Kusner}
\address{Mathematics Department, University of Massachusetts,
      Amherst MA 01003, USA}
\email{kusner@math.umass.edu}

\author[Sullivan]{John M. Sullivan}
\address{Mathematics Department, University of Illinois, Urbana IL 61801, USA}
\email{jms@uiuc.edu}

\begin{abstract}
  We construct the entire three-parameter family
  of embedded constant mean curvature 
  surfaces with three ends and genus zero.
  They are classified by triples of points on the sphere
  whose distances are the necksizes of the three ends.
  Because our surfaces are transcendental, and are not described by
  any ordinary differential equation, it is remarkable to obtain
  such an explicit determination of their moduli space.
  (AMS Classification 2000: 53A10, 58D10)
\end{abstract}
\maketitle

\section*{Introduction}

Surfaces which minimize area under a volume constraint have constant mean
curvature~$H$.  This condition can be expressed as a nonlinear elliptic
partial differential equation:
in terms of a local conformal immersion $f(x,y)$ of the surface into $\R^3$,
we have
$$f_{xx} + f_{yy} = 2 H \, f_x \times f_y,$$
one of the most natural differential equations in geometry.

The case of minimal surfaces, with $H\equiv 0$, reduces
to the Laplace equation.  Thus the coordinates
of a minimal surface are the real parts of holomorphic functions,
and so methods from complex analysis have been exploited to construct
and analyze minimal surfaces, via the Weierstrass representation.

We are interested in surfaces of nonzero constant mean curvature.
For us a \emph{\CMC/ surface} will mean a complete,
properly immersed surface in~$\R^3$, rescaled to have $H\ident1$.
The analysis and construction of \CMC/ surfaces is more delicate than
that of minimal surfaces, since no obvious holomorphic methods are available.

The reflection technique of Alexandrov~\cite{al1} shows that
the unit sphere~$\S^2$ is the only compact embedded \CMC/ surface.
Thus the simplest nontrivial case to consider is embedded surfaces
of \emph{finite topology}.
These are homeomorphic to a compact surface~$\Sigma$ of
genus~$g$ with a finite number~$k$ of points removed.  A neighborhood of each
of these punctures is called an \emph{end} of the surface.

The \emph{unduloids}, embedded \CMC/ surfaces of revolution described in 1841 
by Delaunay~\cite{del}, are genus-zero examples with two ends,
now known to be the only two-ended examples.
The periodic generating curve for each unduloid solves an
ordinary differential equation.
Up to rigid motion, the entire family is parametrized by the
length of the shortest closed geodesics: this \emph{necksize}
ranges from $0$ at the singular chain of spheres, to~$\pi$ at the cylinder.

Unlike the case of minimal surfaces,
continuous families of embedded \CMC/ surfaces may develop self-intersections
not only at the ends, but also on compact sets.
Therefore, to study \CMC/ moduli spaces,
it is natural to consider a more general class of
surfaces, as introduced by Alexandrov~\cite{al2}.
A compact immersed surface is \emph{\alex-emb/}
if it bounds a compact immersed three-manifold.
This is the class of compact surfaces to which
the Alexandrov reflection technique applies,
but the technique can also apply to certain noncompact surfaces.
\begin{definition}
  A \CMC/ surface $M$ of finite topology is \emph{\alex-emb/}
  if $M$ is properly immersed, if each end of $M$ is embedded, and if
  there exists a compact three-manifold~$W$ with boundary $\d W =: \Sigma$
  and a proper immersion $F\takes W\setm\{q_1,\ldots,q_k\}\to\R^3$ whose
  boundary restriction $f\takes \Sigma\setm\{q_1,\ldots,q_k\} \to\R^3$
  parametrizes $M$.
  Moreover, we require that the mean-curvature normal of $M$ points into $W$.
\end{definition}

We define a \emph{triunduloid} to be an \alex-emb/ \CMC/ surface with
genus zero and three ends.
Previously, only very special triunduloids were known to exist:
Kapouleas~\cite{kap} and Mazzeo and Pacard~\cite{mp} used singular
perturbation methods to obtain certain triunduloids with very small necksizes,
while Grosse-Brauckmann~\cite{kgb} constructed triunduloids with high symmetry.
Triunduloids are the basic building blocks for \alex-emb/ \CMC/ 
surfaces with any number of ends.  

Our main result here is the classification 
and construction of all triunduloids.
In contrast, most other classification results for embedded
\CMC/ or minimal surfaces have amounted to proving the uniqueness
of previously known surfaces.
For \alex-emb/ \CMC/ surfaces,
in addition to the round sphere being the only compact example~\cite{al1},
there are none with a single end~\cite{mee},
and only the unduloids have two ends~\cite{kks}.
The situation is similar for embedded minimal surfaces with at least two ends.
(For these, Collin~\cite{col} showed that finite topology is equivalent
to finite total curvature.)
When such surfaces have two ends~\cite{sch} or genus zero~\cite{lr} 
they must be catenoids, while the surfaces with genus one and three ends 
form a one-parameter family~\cite{cos,hk}.

To understand our classification of triunduloids,
note first that a triunduloid always has a mirror symmetry,
unlike general \CMC/ surfaces with more than three ends.
By a theorem of Meeks~\cite{mee}, any triunduloid is contained in a half-space.
Korevaar, Kusner and Solomon~\cite{kks} extended Alexandrov's reflection
technique to \alex-emb/ \CMC/ surfaces contained in a half-space.
It follows that any such surface $M$ is \emph{Alexandrov-symmetric}
in the following sense:
$M$ has a reflection plane~$P$ such that $M\setm P$ consists of two connected
components $M^{\pm}$ whose Gauss images~$\nu(M^{\pm})$ are
contained in the open hemispheres~$\S^2_{\pm}$; moreover $\nu(M\cap P)$
is contained in the equator.
Together with the asymptotics result
of~\cite{kks} this gives for triunduloids (of any genus):
\begin{theorem}[\cite{kks}]\label{thkks} 
  A triunduloid is Alexandrov-symmetric.
  Each of its ends is exponentially asymptotic to a Delaunay unduloid,
  and thus has an asymptotic axis line and necksize $n\in(0,\pi]$.
\end{theorem}

Our Main Theorem below classifies all triunduloids, up to rigid motion.
To avoid orbifold points in the moduli space, we label
the ends of each triunduloid before passing to the quotient space:
\begin{definition}
    The \emph{moduli space $\M$ of triunduloids} consists of all proper,
    \alex-emb/ \CMC/ immersions of $\S^2\setm \{q_1,q_2,q_3\}$ in $\R^3$,
    modulo diffeomorphisms of the domain fixing each $q_i$,
    and modulo orientation-preserving isometries of $\R^3$.
    We define the topology on $\M$ as follows: A sequence in $\M$
    converges if there are representative surfaces
    which converge in Hausdorff distance on every compact subset of $\R^3$.
\end{definition}
The classifying space for triunduloids is quite explicit:
\begin{definition}
    The space of all ordered \emph{triples}
    of distinct points in $\S^2$, up to rotation, is denoted by
    $$\T:=\{(p_1,p_2,p_3)\in \S^2\times\S^2\times\S^2 \colon p_1\not=p_2\not=p_3\not=p_1\}
    \,\big/\, \SO(3).$$
\end{definition}
\noindent
We observe that $\T$ is an open three-ball:
Simply rotate a triple until it lies on a circle of latitude, with $p_1$ on a
fixed longitude, followed by~$p_2$ then $p_3$ proceeding eastward;
the common latitude and the longitudes of $p_2$ and~$p_3$ give coordinates.
Alternatively, we may view $\T$ as the group of conformal motions
modulo rotations, i.e., as hyperbolic space.

It is known~\cite{kmp} that the moduli space~$\M$ of triunduloids is
a real-analytic variety.   We show $\M$ is in fact a manifold,  
by establishing a homeomorphism between the spaces $\M$ and~$\T$:
\begin{maintheorem}
  There is a homeomorphism $\Psi$ from the moduli space~$\M$ of triunduloids
  onto the open three-ball~$\,\T$ of triples in the sphere~$S^2$.
  The asymptotic necksizes of a triunduloid $M\in\M$ are
  the spherical distances of the triple~$\Psi(M)\in\T$.
\end{maintheorem}
\noindent
Thus spherical trigonometry can be used to give necessary conditions
on triunduloids:
\begin{cornn}
  Three numbers $n_1$, $n_2$, $n_3$ in $(0,\pi]$ can be
  the asymptotic necksizes of a triunduloid
  if and only if they satisfy the spherical triangle inequalities
  \begin{equation*} \begin{split}
    n_1 \le n_2+n_3, & \quad
    n_2 \le n_3+n_1, \quad
    n_3 \le n_1+n_2, \\
    & n_1+n_2+n_3 \le 2\pi.
  \end{split} \end{equation*}
  In particular, at most one end of a triunduloid can be cylindrical,
  that is, have necksize~$\pi$.
\end{cornn}
\noindent
When all the inequalities are strict, the triple in~$\T$ has a distinct
mirror image; thus there are exactly two triunduloids with the
same three necksizes (see \figr{triunds}).
\figsix{Three triunduloids and their classifying triples.}
{Three typical triunduloids (top),
and their classifying spherical triples (bottom).
The triunduloid in the middle column
has necksizes $\tfrac\pi2$, $\tfrac{2\pi}3$ and $\tfrac{5\pi}6$; since these
sum to the maximal $2\pi$, its classifying triple is equatorial.
In all three columns, the points of the triples
lie on these same meridians, but on the left and right,
the common latitude is $\pm\tfrac\pi4$ instead of $0$.
These two mirror-image triples correspond to a pair of triunduloids
with a common set of three necksizes,
the smallest of which is $\tfrac\pi3$.
}
The necksize bounds, in the special case $n_1=n_2$,
were described earlier in~\cite{gk}.

The forces of the ends (as defined in Section~\ref{se:prop})
of any triunduloid must balance.
This leads to explicit trigonometric formulas~\cite{gks1} 
which prescribe the angles between the asymptotic axes of the ends
in terms of the necksizes.
(We lack, however, equally good control over one additional
parameter describing the asymptotics of each end,
namely the phase or position of the asymptotic necks along the axis.)
Force balancing also gives another intuitive explanation for
why $\M$ is three-dimensional:
The three force vectors give nine parameters, while force balancing 
and the action of~$\SO(3)$ each reduce the count by three.

\subsection*{Outline of the proof}
The main steps in the proof of our theorem (see also~\cite{pnas})
are to define the classifying map $\Psi$ from triunduloids to spherical
triples, and then to prove it is injective, proper, and surjective.

In Section~\ref{se:class} we define $\Psi$.
We make use of Alexandrov symmetry to decompose a triunduloid in~$\M$
into two halves.
Each half is simply connected, so Lawson's conjugate-cousin
construction~\cite{law} gives a minimal cousin in the three-sphere~$\S^3$.
The three boundary components of the cousin are fibers in a single
Hopf fibration~\cite{kar},
and hence they project under the Hopf map to the three points of a
spherical triple in~$\T$.

Our approach to conjugate cousins, explained in Section~\ref{se:coco},
is different from Lawson's.  It allows us to show that the minimal cousin itself
Hopf-projects to an immersed disk in~$\S^2$; the spherical metric
induced on this disk depends only on the triple in~$\T$.
Pulling back the Hopf bundle, we get a circle-bundle over this disk,
and we regard the minimal cousin as a section.
The cousin of any other triunduloid with the same triple
gives another minimal section of the same bundle,
and we use the maximum principle to show these sections are congruent.  
This proves the injectivity of~$\Psi$, the main result of Section~\ref{se:inj}.

To show $\Psi\colon\M\to\T$ is a homeomorphism, we employ a 
novel continuity method, applicable to problems with real analytic
spaces of solutions.
Instead of showing $\Psi(\M)$ is closed as a subset of~$\T$,
we show that our injective map $\Psi$ is continous and proper,
using the area estimate from~\cite{kk}
and extending their curvature estimate (Section~\ref{se:prop}).

Instead of showing $\Psi(\M)$ is open we 
use an interesting topological argument 
to prove $\Psi$ is surjective and a homeomorphism (Section~\ref{se:surj}).
For this, we depend on the fact that our moduli space~$\M$ of triunduloids
is locally a real analytic variety~\cite{kmp}.
A real-analytic variety is known to have only even-fold branching.  
But $\Psi$ maps the variety~$\M$ injectively to the manifold~$\T$ and so
the branching can only be two-fold.  Hence $\M$ is a topological manifold.
Since we can show $\M$ has dimension three (using~\cite{kmp,mr,mp}), 
while $\T$ is connected, our injective map is surjective.

\subsection*{Open problems and related work}
A triunduloid is called \emph{nondegenerate} when 
all its $L^2$~Jacobi fields vanish (see Section~\ref{se:surj}).  From
general theory~\cite{kmp} it is known that near a nondegenerate 
triunduloid, the moduli space locally must be a three-manifold.
Our classification theorem suggests that all triunduloids
may be nondegenerate. 
It would be interesting to prove this conjecture,
since Ratzkin has recently shown~\cite{ratz}
that nondegenerate triunduloids can be glued end-to-end
to construct genus-zero \CMC/ surfaces with any finite number of ends.
Presumably, \CMC/ surfaces of any finite genus~$g$
can be constructed similarly by gluing triunduloids.

More generally one can consider \alex-emb/ \CMC/ surfaces with $k$~ends 
whose asymptotic axes lie in a common plane.  Our methods generalize to
show that the moduli space of these
\emph{coplanar $k$-unduloids} of genus~$0$ is also a manifold
(of dimension $2k-3$), though in general not contractible.
Because this case is technically more involved, we devote a separate
paper~\cite{kund} to it.

Cos{\'\dotlessi}n and Ros~\cite{cr} adapted our approach
to classify coplanar minimal $k$-noids in~$\R^3$.
Again, they use boundary data of the surface conjugate to the upper half
as a classifying map.

Surfaces of mean curvature one in hyperbolic space are more
similar to minimal surfaces in Euclidean space than to \CMC/ surfaces,
thanks to the cousin relationship exploited by Bryant~\cite{bry}.
This allowed Umehara and Yamada~\cite{uy} to prove a result like ours for
mean-curvature-1 trinoids in hyperbolic space.

The DPW representation~\cite{dpw} uses loop group techniques to parameterize
\CMC/ surfaces.
In work in progress, Dorfmeister and Wu as well as Schmitt (see~\cite{kms})
have used this representation to exhibit genus-zero \CMC/ surfaces
with three ends.  It seems hard to determine, however, whether any of
these are Alexandrov-embedded.

\subsection*{Acknowledgements}
We wish to thank Claus Hertling for his advice on real analytic
varieties, Frank Pacard for many useful suggestions,
and Hermann Karcher for helpful discussions.

We also gratefully acknowledge the hospitality of
the Institute for Advanced Study
and the Center for GANG (through NSF grant DMS~96-26804)
at various times during this project, as well as that
of the CMLA at \'Ecole Normale Sup\'erieure Cachan (K.G.B.)
and the NCTS in Hsinchu, Taiwan (R.B.K.).
This work was partly supported by the Deutsche Forschungsgemeinschaft
through SFB 256 (K.G.B.), and by the National Science Foundation
under grants DMS~97-04949 and DMS~00-76085 (R.B.K.)
and DMS~97-27859 and DMS~00-71520 (J.M.S.).

\section{Conjugate cousins}\label{se:coco}
\subsection{Conjugate minimal surfaces}

A simply connected minimal surface $M$ in $\R^3$ has,
as is well known, a conjugate minimal surface $\Mconj$.
These surfaces are isometric, and have a common Gauss map~$\nu$:
at corresponding points they have the same tangent plane.
However, tangent vectors are rotated by $90^\circ$ under the isometry.

In terms of a parametrization $f\colon\Omega\to\R^3$ of $M$
(where $\Omega\subset\R^2$) and the corresponding parametrization
$\fconj$ of~$\Mconj$, conjugation can be described as follows.
Since $f$ and $\fconj$ are isometric, they induce the same pullback
metric on $\Omega$.  We let $J$ denote rotation by $90^\circ$
with respect to this metric.
Then the parametrizations are related by the first-order
differential equation
\begin{equation}\label{eq:first}
	\dfconj= df\circ J.
\end{equation}
Equivalent formulations of this equation include
$\dfconj=-\star df$, using the Hodge star induced on $\Omega$,
and $\dfconj = \nu\times df$,
using the cross product in $\R^3$ with the unit normal vector $\nu$.

If $f$ is conformal, $J$ is the ordinary rotation in $\R^2$, and
\eqn{first} becomes the Cauchy-Riemann system
$$
  \fconj_x = f_y, \qquad \fconj_y = -f_x.
$$
That is, when a minimal surface $M$ is given a conformal parametrization,
its coordinate functions are harmonic, and those of $\Mconj$ are the
conjugate harmonic functions.

Symmetry curves on conjugate surfaces have a striking relationship.
If $f\circ\gamma$ is a curve in $M$ (parametrized by arclength)
then~\eqn{first} implies that its
\emph{conormal} vector $\eta:=df(J\gamma')$ equals the tangent
vector $\dfconj(\gamma')$ of the conjugate curve in~$\Mconj$.
Thus $f\circ\gamma$ has constant conormal $u\in\S^2$ if and only if
the curve $\fconj\circ\gamma$ is contained in a straight line in
the direction $u$;
that is, an \emph{arc of planar reflection} in~$M$ corresponds to an
\emph{arc of rotational symmetry} in~$\Mconj$.

According to Bonnet's fundamental theorem,
a surface in space is determined up to rigid motion by second-order
data: its metric and its shape operator.
This leads to an equivalent, second-order description of conjugation:
$M$ and~$\Mconj$ are isometric, and their shape operators are related by
$J\circ\Sconj = S$.
The Gauss and Codazzi equations, which are the integrability conditions in
Bonnet's theorem, confirm that the conjugate surface $\Mconj$ exists exactly
when $M$ is minimal.

\subsection{Conjugate cousins and quaternions}

Lawson~\cite{law} described a similar notion of \emph{conjugate cousin}
surfaces, where $M$ is a \CMC/ surface in~$\R^3$ and $\Mconj$
is an isometric minimal surface in~$\S^3$.
He gave a second-order description, in which the shape operators are
related by $J\circ\Sconj=S-\id$.
Again, the integrability conditions in Bonnet's theorem
prove existence of~$\Mconj$ exactly when $M$ is \CMC/ and, conversely, existence
of~$M$ exactly when $\Mconj$ is minimal~\cite{law,kar,kgb}.
It is the Gauss equation which shifts the unit mean curvature of
$M\subset\R^3$ to ambient sectional curvature of $\S^3\supset\Mconj$.
(With the terminology \emph{cousin} we follow Bryant~\cite{bry}, who similarly
relates \CMC/ surfaces in hyperbolic space to minimal surfaces in $\R^3$.)

Here we will give a description of conjugate cousins via
a first-order system of differential equations.
When properly interpreted, it says again that
the surfaces share a common normal,
and that tangent vectors rotate by a right angle.
This description seems to have been mentioned first in~\cite{gbp}
(see also~\cite[p.66]{hel}),
while the first proof (in a different context) appears in~\cite{obe}.
Karcher's earlier work~\cite{kar} can be regarded as a first-order
description of the symmetry curves alone.

System~\eqn{first} for minimal surfaces in~$\R^3$ equates tangent vectors at
two different points of~$\R^3$, implicitly making use of parallel translation.
In the case of conjugate cousins, to compare tangent vectors to~$\R^3$
and~$\S^3$, we use the Lie group structure of~$\S^3$.
If $L_p$ denotes left translation of $\S^3$ by $p\in\S^3$, then
we use the left translation of tangent vectors,
$dL_p\colon T_1\S^3\to T_p\S^3$, to generalize~\eqn{first} to
\begin{equation}\label{eq:cousin}
	\dfconj=dL_\fconj \, (df\circ J).
\end{equation}
That is, given a tangent vector to the \CMC/ surface $f\takes \Omega\to\R^3$,
we first rotate it by $90^\circ$ within the tangent plane.
This rotated vector is an element of $T_{f(z)}\R^3 \isom \R^3$,
which we identify with $T_1\S^3$.  Finally,
we left-translate by~$dL_{\fconj(z)}$ to obtain a vector in~$T_{\fconj(z)}\S^3$.

It is convenient to view the Lie group~$\S^3$ as the set of unit quaternions.
Then tangent vectors to~$\S^3$ are also represented by quaternions,
and left translation (of points or tangent vectors) is simply multiplication.
We identify $\R^3$ with $T_1\S^3=\Im\H$, the imaginary quaternions.
This way, \eqn{cousin} becomes the quaternion-valued equation
$\dfconj=\fconj\,df\after J$.
Again, if $f$ (or equivalently $\fconj$) is conformal, we can write this as
\begin{equation}\label{eq:fstHcnf}
  \fconj_x = \fconj\,f_y, \qquad \fconj_y = -\fconj\,f_x.
\end{equation}
In order to investigate the integrability of this system, let us now
briefly review some of the geometry of surfaces in $\R^3=\Im\H$ and in
$\S^3\subset\R^4=\H$ using quaternionic notation.

For imaginary quaternions, $p,q\in\Im\H$, we can express the quaternionic
product in terms of the inner product and cross product:
$\Re(pq)=-\langle p, q\rangle\in\Re\H=\R$
and $\Im(pq)= p\times q\in\Im\H=\R^3$, giving
$pq=-\langle p, q\rangle + p\times q$.
In particular, imaginary $p$ and $q$ commute if and only if they are
parallel, and anticommute if and only if they are orthogonal.
Therefore, conformality for a map $f\takes\Omega\to\Im\H$ can be expressed as
\begin{equation}\label{eq:conf}
  |f_x|^2=|f_y|^2, \qquad f_xf_y=-f_yf_x.
\end{equation}

Given a surface in an oriented three-manifold,
a choice of unit normal $\nu$ is equivalent to
a choice of orientation for the surface.
We always orient a \CMC/ surface by choosing the
(inward) \emph{mean curvature normal}.
Once we orient a surface, we will always use oriented parametrizations
such that $f_x\times f_y$ is a positive multiple of~$\nu$.
The $90^\circ$ rotation~$J$ on the domain is then oriented in
the usual counterclockwise manner, so that
$\big(df(X),df(JX),\nu\big)$ is positively oriented for any $X\in T\Omega$.

With the above conventions, the mean curvature~$H$
of a conformally parametrized surface $f\colon\Omega\to\R^3$
can be measured by $\Delta f = 2 H\, f_x\times f_y$, where
$\Delta=\tfrac{\d^2}{\d x^2}+\tfrac{\d^2}{\d y^2}$
is the usual Laplacian on~$\Omega$.
Expressing the cross product quaternionically, we see that
$f$ is \CMC/ (meaning that $H\ident1$) if and only if
\begin{equation}\label{eq:cmc}
   \Delta f = 2f_xf_y .
\end{equation}
Similarly, a conformal $\fconj\takes\Omega\to\S^3$ is minimal
if and only if it is harmonic:
\begin{equation}\label{eq:minimal}
   \Delta \fconj = -\fconj\,|\dfconj|^2,
\end{equation}
where $|\dfconj|^2 = |\fconj_x|^2 + |\fconj_y|^2$.

\subsection{Existence of cousins}

We are now ready to prove the theorem giving the existence of
isometric conjugate cousin surfaces $M\subset\R^3$ and $\Mconj\subset\S^3$.
We will find the curvature equation~\eqn{cmc} for~$f$ is the integrability
condition for~$\fconj$, and similarly~\eqn{minimal} is the one for~$f$.
Thus cousins exist exactly when $f$ is \CMC/ or $\fconj$ is minimal.

\begin{theorem}\label{thcous}
  Let $\Omega$ be a simply connected domain in $\R^2$.
  Then, for each immersion $f\takes\Omega\to\R^3=\Im\H$ of
  constant mean curvature~$1$,
  there is an isometric minimal immersion $\fconj\takes\Omega\to\S^3$,
  determined uniquely up to left translation, satisfying
  \begin{equation}\label{eq:firstH}
    \dfconj = \fconj \,df\circ J.
  \end{equation}
  Conversely, given an oriented minimal~$\fconj$, there is a
  \CMC/ immersion~$f$, unique up to translation, satisfying~\eqn{firstH}.
  In either case, the surface normals agree, in the sense that
  $\nuconj=\fconj \nu$.
\end{theorem}
\begin{proof}
  The results are invariant under
  orientation-preserving diffeomorphisms of~$\Omega$.
  Thus it is sufficient to prove the theorem
  assuming a conformal parametrization.

  Because \eqn{firstH} implies that tangent vectors to the two surfaces
  have the same lengths, and lie within $T\S^3$ or $T\R^3$
  as appropriate, it is clear that solutions are isometric surfaces
  in $\S^3$ and $\R^3$.
  Equivalently, it is straightforward to check algebraically
  that any quaternionic solutions to \eqn{firstH} have
  the following two properties:
  first, $|\fconj|$ is constant if and only if $\Re f$ is constant;
  second, the metrics induced by $f$ and $\fconj$ are conformal,
  with lengths scaling by the factor $|\fconj|$.

  Now assume $\fconj$ exists and write \eqn{firstH} as
  $df=-\fconj^{-1}\dfconj\!\circ\! J$.
  This is integrable for~$f$ exactly when
  $$
    d(-\fconj^{-1}\dfconj\!\circ\! J)
    = -d(\fconj^{-1}) \wedge \dfconj\!\circ\! J
      - \fconj^{-1}d(\dfconj\!\circ\! J)
  $$
  vanishes.  In conformal coordinates, the first term becomes
  $$
    \fconj^{-1}\dfconj \,\fconj^{-1} \wedge \dfconj\!\circ\! J
     = \fconj^{-1} (\fconj_xdx+\fconj_y dy)
       \wedge \fconj^{-1} (\fconj_ydx-\fconj_x dy)
    = -\big((\fconj^{-1}\fconj_x)^2 +(\fconj^{-1}\fconj_y)^2\big)
       dx\wedge dy
  $$

  Now $\fconj^{-1}\fconj_x$ and $\fconj^{-1}\fconj_y$ are pure imaginary,
  being tangent to $\S^3$ at $1$.  It follows that
  $$
    -\big(\fconj^{-1}\fconj_x\big)^2- \big(\fconj^{-1}\fconj_y\big)^2
    = |\fconj^{-1}|^2\big(|\fconj_x|^2+|\fconj_y|^2\big)
    = |\dfconj|^2
  $$
  and so 
  $$
    d\big(-\fconj^{-1}\dfconj\circ J\big)
    = \big(|\dfconj|^2  + \fconj^{-1}\Delta\fconj\big)dx\wedge dy.
  $$
  Comparing with~\eqn{minimal}, we see that our system is integrable for $f$
  exactly when $\fconj$ is minimal.

  Conversely, let $f\takes\Omega\to\R^3=\Im\H$ be conformal,
  and set $\alpha=df\after J=f_y\,dx-f_x\,dy$, so that~\eqn{firstH}
  becomes $\dfconj=\fconj\alpha$, or $\fconj^{-1}\,\dfconj=\alpha$.
  Using the machinery of Maurer--Cartan, the
  integrability condition for such an equation is
  $$
    d\alpha + \alpha\wedge\alpha = 0.
  $$
  (See~\cite[Chap.~3]{sharpe} for a discussion
  of this fundamental theorem of calculus for Darboux derivatives.)
  But we calculate $d\alpha=-\Delta\! f\,dx\wedge dy$, and
  $\alpha\wedge\alpha = f_xf_y\,dx\wedge dy + f_yf_x\,dy\wedge dx
  = 2f_xf_y\,dx\wedge dy$, using~\eqn{conf}.
  Thus the integrability condition is precisely~\eqn{cmc},
  so a solution $\fconj$ exists exactly when $f$ is \CMC/,
  and it is determined up to left translation.

  In either case, $f$ and $\fconj$ are conjugate cousins; either
  could be derived from the other via~\eqn{firstH},
  and both integrability conditions are satisfied.
  This verifies that $\fconj$ is minimal and $f$ is \CMC/,
  no matter which we started with.
\end{proof}

The simplest example of cousins is the unit sphere~$\S^2$ contained
both in $\S^3$ and $\Im\H$ (as their intersection).
It is its own conjugate cousin via the identity map.
Indeed, if $f$ parametrizes $\S^2$ then
$f\,df\circ J=f\times df\circ J=df$ and so $\fconj=f$.
In this case, the point $1\in\S^3$ and its left translate
$f\in\S^2\subset\S^3$ always have spherical distance~$\tfrac\pi 2$;
therefore, on the level of tangent vectors a $90^\circ$ rotation is induced
which cancels the effect of~$J$.

For the record, we show that our first-order description of conjugate
cousins agrees with Lawson's original second-order description~\cite{law}:
\begin{proposition}
  The shape operators $S$ of~$f$ and $\Sconj$ of~$\fconj$
  are related by $J\after\Sconj=S-\id$.
\end{proposition}
\begin{proof}
  By definition of the shape operator,
  $\dfconj\after \Sconj=-d\nuconj = -\dfconj\,\nu - \fconj\, d\nu$,
  using $\nuconj=\fconj\nu$.
  Substituting $\dfconj=\fconj\,df\after J$ and $d\nu=-df\after S$, we get
  $$
    df\after J \after \Sconj(\cdot) =
    -\big(df\after J(\cdot) \big)\,\nu + \,df\after S(\cdot).
  $$
  The first term on the right can be simplified:
  $df\after J$ takes values in the tangent plane to $M$;
  since multiplication with the normal $\nu$ simply rotates
  by $90^\circ$, we get
  $$\big(df\after J(\cdot)\big)\nu = \big(df\after J(\cdot)\big)\times\nu
   = df(\cdot).$$
  Therefore, $J\after\Sconj=-\id + S$ as desired.
\end{proof}

\subsection{Hopf fields and spherical cousins}

Left-translating a unit vector $u\in\S^2\subset\Im\H=T_1\S^3$
gives us a left-invariant
vector field $p\mapsto pu$ on $\S^3$, called the \emph{$u$-Hopf field}.
Its integral curve through $p\in\S^3$ is a \emph{$u$-Hopf circle},
the great circle through~$p$ and~$pu$ parametrized by
$pe^{tu} := p(\cos t+u\sin t)$.
The \emph{$u$-Hopf projection} $\Pi_u\takes\S^3\to\S^2$ is
defined by $\Pi_u(p) = pu\bar p$.  This projection is a fibration;
its fibers are the $u$-Hopf circles, since
$$\Pi_u\big(p(\cos t+u\sin t)\big)
= p(\cos t+u\sin t)u(\cos t-u\sin t)\bar p
= \Pi_u(p).$$

Later, we will also consider the $u$-Hopf projection of
a $v$-Hopf circle~$\gamma$ with $\langle u,v \rangle =\cos\theta$.
It is a round circle $\Pi_u(\gamma)$ of spherical
radius~$\theta$, covered twice at constant speed $2\sin\theta$.
Indeed, the $u$-Hopf circles through the different points of $\gamma$
foliate a distance torus~$T$ in~$\S^3$.
They are curves of homology class $(1,1)$ on~$T$, while $\gamma$ is $(1,-1)$.
The two polar core circles of~$T$ are $u$-Hopf, and project under~$\Pi_u$
to the two antipodal centers of $\Pi_u(\gamma)$ in~$\S^2$.

As on conjugate minimal surfaces,
symmetry curves on cousin surfaces are related.
This time, using \eqn{firstH}, we find that
a curve $f\circ \gamma$ in~$M$ with constant conormal $u\in\S^2$
corresponds to a curve $\fconj\circ \gamma$ in~$\Mconj$
with tangent vector $\fconj u$ in the $u$-Hopf field:
\begin{proposition}\label{pr:bdycurves}
  Let $M$ and $\Mconj$ be conjugate cousins.
  A curve $f\circ\gamma$ in~$M$ is a curve of planar reflection for~$M$,
  contained in a plane with unit normal $u\in\S^2$,
  if and only if the curve $\fconj\circ\gamma$ in $\Mconj$ is contained
  in a $u$-Hopf circle in~$\S^3$.
\end{proposition}

Since left translation is an isometry for
the metric on $\S^3$ we conclude from $\nuconj=\fconj\nu$ that
$\langle \nuconj,\fconj u\rangle_{\S^3} = \langle \nu,u \rangle_{\R^3}$
on~$\Omega$.
When this is nonzero, we get:
\begin{proposition}\label{pr:trvs}
  A \CMC/ surface $M$ is transverse to the foliation of~$\R^3$
  by lines in the $u$~direction if and only if its minimal cousin
  $\Mconj$ is transverse to the foliation of~$\S^3$ by $u$-Hopf circles.
  In this case, the Hopf projection $\Pi_u$ immerses $\Mconj$ to $\S^2$.
\end{proposition}

\subsection{Unduloids and their helicoid cousins}

The \emph{spherical helicoid} is parametrized by $h\colon\R^2\to\S^3$,
\begin{align*}
  h(x,y)
         :=& \cos x\,\cos ny - \i\,\cos x\,\sin ny
            + \j\,\sin x\,\cos(2\pi-n)y + \k\,\sin x\,\sin(2\pi-n)y,
\end{align*}
where $\{1, \i, \j, \k\}$ is the standard basis of $\H=\R^4$
and $n$ is a real parameter.
For $n\not=0,2\pi$ the helicoid $h$ is an immersion;
moreover it is minimal (see \cite[\paragraph 2]{kgb}
and \cite[\paragraph 7]{law}).  

To analyse the symmetries of the cousin, we rewrite $h$ as
\begin{equation}\label{eq:helic}
  h(x,y)
  = p(y)\big(\cos x+u(y)\sin x\big),
\end{equation}
where $p(y) := \cos ny-\i\,\sin ny   \in\S^3$
and $u(y) := \j\,\cos 2\pi y + \k\,\sin2\pi y \in\S^2\subset\Im\H$.
The $x$-coordinate lines $x\mapsto h(x,y)$
trace out the $u(y)$-Hopf circles through $p(y)$.
These great-circle rulings are orthogonal to the \emph{axis} of the helicoid,
$y\mapsto h(0,y)=p(y)$, along the $(-\i)$-Hopf circle through~$1$.
In fact, the helicoid has another axis,
$$
  y\mapsto h(\tfrac\pi 2,y)
  = p(y)u(y)
  = \j\,\cos(2\pi-n)y + \k\,\sin(2\pi-n)y
$$
along the $(-\i)$-Hopf circle through~$\j$;
this is the polar circle of the original axis.
Each ruling intersects both axes, and its tangent vector at the first axis
coincides with its position on the polar axis:
$\tfrac{\d h}{\d x}(0,y)=h(\tfrac\pi 2,y)$.
Since the polar axis is parametrized with constant speed,
the tangent vectors rotate with constant speed along the first axis,
thus justifying the name spherical helicoid.

Our parametrization $h(x,y)$ is along orthogonal asymptotic lines,
confirming that~$h$ is a minimal surface.
Indeed, the $x$-coordinate lines are the great-circle rulings.
The $y$-coordinate lines are helices contained in
distance tori around the polar axes of~$h$.
Within the torus, the helix is geodesic and so
its curvature vector in~$\S^3$ is normal to the torus, and hence
tangent to the helicoid~$h$.

The next lemma will show that the cousins of these helicoids
are the \CMC/ surfaces of revolution characterized by Delaunay~\cite{del}:
for $n\in(0,\pi]$, the helicoid~$h$ is the cousin of
an embedded unduloid of necksize~$n$;
for $n<0$, it is the cousin of an immersed \emph{nodoid}.
Indeed, by Proposition~\ref{pr:bdycurves} the cousin \CMC/ surface
is foliated by arcs of planar reflection; the normals to the mirror planes are
the tangent Hopf fields $u(y)$ of the rulings.
Since $u(y)$ rotates through all vectors perpendicular to~$\i$, the cousin is
a surface of revolution around the $\i$-axis.

\begin{lemma}\label{le:unduloids}
  If we fix $n\in(0,\pi]$, the spherical helicoid $h(x,y)$,
  given by~\eqn{helic} for $x\in\R$ and $y\in(-\tfrac 14,\tfrac 14)$,
  is a minimal surface, the cousin $\plustilde U$
  of a half-unduloid $U^+$ with necksize~$n$, whose boundary $\d U^+$
  is contained in the $\i\j$-plane.
  Moreover, the Hopf projection $\Pik$ immerses $\plustilde U$ into the
  two-sphere $\S^2$ in such a way that the boundary
  $\Pik(\d\plustilde U)$
  consists of two points at spherical distance~$n$.
\end{lemma}
\begin{proof}
  The cousin of $h$ is a \CMC/ surface of revolution with meridians contained
  in mirror planes with normal $u(y)$.
  Since $y\in(-\tfrac14,\tfrac14)$ is a maximal interval
  for which these planes are nonhorizontal, the upper half $U^+$ is
  indeed conjugate to~$h$ on the claimed domain.

  The axes of the helicoid give circles of planar reflection
  on the \CMC/ surface, contained in planes with normal~$\i$.
  In particular, the first axis is parametrized with speed $n$
  and so gives a neck circle with circumference~$n$, while the polar axis
  is parametrized with speed $2\pi-n$ and so gives a circle with larger
  circumference $2\pi-n$.
  Since $n>0$, this latter bulge circumference is less than $2\pi$, showing
  that the \CMC/ surface is an unduloid, not a nodoid.
  Note that this description shows easily that the meridians of any unduloid
  have length $\tfrac\pi2$ between neck and bulge.

  Since unduloids are embedded, $U^+$ is transverse to the vertical lines,
  and so $\Pik$ immerses $\plustilde U$ by Proposition~\ref{pr:trvs}.
  The bounding rulings
  $h(x,\pm\tfrac14) =p(\pm\tfrac14)(\cos x \pm \k\,\sin x)$
  are $(\pm\k)$-Hopf circles through $h(0,\pm\tfrac14)=p(\pm \tfrac14)$,
  and hence \kHopf/ project to single points
  $\k\,\cos\tfrac n2 \mp \j\,\sin\tfrac n2$,
  which are at spherical distance~$n$.
\end{proof}

\section{The classifying map}\label{se:class}

Our triunduloid moduli space $\M$ consists of
equivalence classes of triunduloids up to rigid motion.
For convenience, we assume the triunduloid $M$ representing a class
has been rotated as follows: The symmetry plane $P$,
guaranteed by Theorem~\ref{thkks}, is the $\i\j$-plane in~$\R^3$;
moreover, the asymptotic axis of the first end points in the $\i$-direction,
and successively labeled ends occur in counterclockwise order in~$P$.
For each class in~$\M$, this determines the representative~$M$
uniquely up to horizontal translation.
The intersection of $M$ with the open half-space above the symmetry plane
will be denoted $M^+$, the \emph{upper half} of $M$.

Since $M$ has genus zero, $M^+$ is simply connected.
Hence Theorem~\ref{thcous} gives a conjugate cousin $\plustilde M\subset\S^3$,
which is determined uniquely up to left translation.
Since $\d M^+$ consists of three arcs of reflection, all contained in
the $\i\j$-plane~$P$ with normal $\k\in\S^2$,
it follows from Proposition~\ref{pr:bdycurves} that
$\plustilde M$ is bounded by three arcs of \kHopf/ circles.
Thus $\Pik(\d\plustilde M)$ consists of three points in~$\S^2$.

We want to use this triple to define our map $\Psi\colon\M\to\T$,
so we need to check that the triple is well defined up to rotation,
and that the three points are distinct.
But~$\plustilde M$, and hence its boundary, is well defined up
to left translation by some $a\in \S^3$,
and $\Pik(ap)= ap\,\k\,\bar p\,\bar a = a\Pik(p)\bar a$.
This is indeed a rotation of~$\S^2$ applied to $\Pik(p)$.
We show the triple is distinct by proving the second part of our Main Theorem.
\begin{theorem}\label{th:necksizes}
  The spherical distances between the points of the triple
  $\Pik(\d\plustilde M)$ are the necksizes of the triunduloid~$M$.
\end{theorem}
\begin{proof}
  Let $E_i$ be the $i^{\text{th}}$ end of $M$, and $E_i^+=E_i\cap M^+$.
  Theorem~\ref{thkks} states that $E_i$ is exponentially
  asymptotic to a Delaunay unduloid~$U$ of necksize~$n_i$.
  By \eqn{firstH}, conjugation preserves exponential convergence,
  and so the bounding rays of $\plustilde{E_i}$ are asymptotic to the
  great circles bounding the helicoid~$\plustilde U$.
  However, by Proposition~\ref{pr:bdycurves} the two rays
  in the boundary of $\plustilde{E_i}$ are themselves great-circle rays.
  Thus they must cover the great circles bounding the cousin of the asymptotic
  unduloid.
  Lemma~\ref{le:unduloids} shows that these Hopf-project to points
  at spherical distance~$n_i$.
\end{proof}

\begin{definition}
The \emph{classifying map} for triunduloids is given by
$$\Psi\colon\M\to\T, \qquad
  \Psi(M):=\Pik(\d \plustilde M )\in \T.$$
\end{definition}
\noindent
By Theorem~\ref{th:necksizes} and the discussion above, $\Psi$ is
well defined with values in~$\T$.
In Section~\ref{se:prop}, we will use curvature bounds to
check that $\Psi$ is continuous and proper.

\begin{remark}
Equally well, we could define $\Psi$ on a triunduloid $M$ with arbitrary
mirror plane~$P$, by Hopf-projecting in the direction $u$ normal to $P$.
If $M$ is rotated by $A\in\SO(3)$, then its cousin is also rotated by~$A$,
extended to act on $\S^3$ by fixing~$1$ and preserving
$\S^2\subset\R^3=\Im\H$.
Hopf projection is equivariant with respect to such rotations:
if $A$ is realized through conjugation by a unit quaternion $a$, then
$$
  \Pi_{Au}(Ap)
  =(ap\bar a)(au\bar a)\bar{(ap\bar a)}
  =a(pu\bar p)\bar a
  =A\big(\Pi_u(p)\big).
$$
Thus surfaces representing the same class in $\M$ have Hopf images
related by a rotation, and indeed give the same triple in~$\T$.
\end{remark}

\section{Injectivity}\label{se:inj}

The injectivity of our classifying map $\Psi$ is a uniqueness
theorem for triunduloids with the same classifying triple.
By definition of the classifying map, their cousin disks
are bounded by the same three Hopf fibers;
we will apply the maximum principle to show
the cousins are uniquely determined by these bounding fibers.

It is well known that two minimal graphs over the same domain
in $\R^2$ coincide, provided they have the same boundary data.
Indeed, by a vertical translation, the graphs can be assumed to have a
one-sided contact, in which case they agree by the maximum principle.
We can consider the two graphs as minimal sections of the
trivial (vertical) line bundle over the domain;
this suggests how we will generalize the argument
to certain minimal surfaces in $\S^3$.

As in the previous section,
we assume $M$ is rotated to have horizontal Alexandrov symmetry plane.
Since the cousin minimal disk $\plustilde M\subset\S^3$ is then
transverse to the \kHopf/ circles,
we can use the \kHopf/ projection of $\plustilde M$ to define the
domain disk for a generalized graph.

\begin{proposition}\label{pr:phiimmerses}
  Define $\Phi(M):=\Pik(\plustilde M)$. This is an immersed open disk in~$\S^2$,
  well defined up to rotation of~$\S^2$.
\end{proposition}
\begin{proof}
  By Theorem~\ref{thkks} the immersed open disk $M^+$ is transverse to the
  vertical lines in~$\R^3$.
  It follows from Proposition~\ref{pr:trvs} that
  $\Pik$ immerses $\plustilde M$ into~$\S^2$.
  Finally, $\plustilde M$ is well defined up to left translation,
  so as in Section~\ref{se:class}, $\Phi$ is well defined up to rotation.
\end{proof}
\noindent
Note that, by our conventions, this immersion preserves the orientation
given by the inward normals on $M$ and $\S^2$.

The disk $\Phi(M)$ carries a spherical metric which, as we will show in
section~\ref{ss:sphtri} below, depends on the classifying triple alone.
Therefore the pullback of the Hopf bundle to this disk defines a bundle
whose geometry depends on the classifying triple only.
This allows us, in section~\ref{ss:unique}, to apply the maximum principle
in the bundle to prove the injectivity of~$\Psi$.

\subsection{Spherical metrics on half-unduloids}

A metric on the oriented open two-disk~$D$ is called
\emph{spherical} if it is locally isometric to $\S^2$.
We will also consider the metric completion $\hat D$,
and refer to $\hat D \setm D$ as the \emph{completion boundary} of $D$.
A spherical metric defines an oriented \emph{developing map}
$\phi\colon D \to\S^2$, unique up to rotation,
which extends continuously and isometrically to~$\hat D$.

A particular case we consider is metrics whose completion boundary
includes an arc developing onto a great-circle arc in~$\S^2$ of some
length $n\in(0,2\pi)$.
We can \emph{join} two such metrics if the lengths $n$ are the same,
by developing the two disks so that the images in~$\S^2$ extend one
another across the common boundary arc.  The union of the disks,
together with the open arc, is again a disk with a spherical metric.

Suppose we are given a \emph{minimizing great-circle arc} in~$\S^2$,
that is, a closed arc~$e$ of length $n\in(0,\pi]$.
Its complement is a disk with a spherical metric; we
call this a \emph{slit sphere} of \emph{slit length} $n$.
The completion boundary of a slit sphere
is a loop consisting of two arcs of length $n$.
We can join two, or more generally $k\ge 2$,
slit spheres of the same slit length.
Their join is a chain of $k$ slit spheres,
whose developing map has degree $k$ almost everywhere.
The completion boundary is again a loop of two arcs,
each developing onto the original arc~$e$.

A \emph{ray of slit spheres} is obtained by
successively joining infinitely many slit spheres of common slit length.
The developing map of such a ray has infinite degree.
Its completion boundary is a single closed arc of length $n$.
Joining two rays of slit spheres gives a \emph{line of slit spheres},
whose developing map is the universal covering
of~$\S^2\setm\{p,q\}$, where $p$ and $q$ are at distance~$n$.
The completion boundary can be identified with $\{p,q\}$.

\begin{lemma}\label{le:slitray}
If the completion boundary of a spherical metric on the disk
is a closed arc developing onto a minimizing great-circle arc,
then the metric is a ray of slit spheres.
\end{lemma}
\begin{proof}
By joining two copies of the metric, we obtain a spherical metric
on the disk~$D$ whose completion boundary develops
to the two endpoints~$\{p,q\}$ of the original arc.
Any path in $\S^2\setm\{p,q\}$ then lifts to~$D$,
so $D$~is the universal isometric covering space of
$\S^2\setm\{p,q\}$, that is, a line of slit spheres.
The arc across which we joined the two copies of
the original metric develops to~$\overline{pq}$,
and thus decomposes~$D$ into two rays of slit spheres.
\end{proof}

Given an unduloid~$U$ of necksize~$n$, by Lemma~\ref{le:unduloids}
the Hopf projection of the cousin induces a spherical metric
on~$U^+$.
A~\emph{half-bubble} of an unduloid $U$ is the open subset of the
upper half $U^+$ enclosed by the semicircles around two consecutive
necks.  On a cylinder, we let a half-bubble be the region
enclosed by any two semicircles at distance~$\pi$.
\begin{lemma}\label{le:hopfund}
  Let $U^+$ be the upper half of an unduloid of necksize $n$
  (with boundary in a horizontal plane as in Lemma~\ref{le:unduloids}),
  and let $B^+\subset U^+$ be a half-bubble.
  Then $\Pik(\plustilde B)$ is a slit sphere of slit length~$n$.
  Consequently, $\Pik(\plustilde U)$ is a line of slit spheres.
\end{lemma}
\begin{proof}
  The unduloid $U$ is the cousin of the
  spherical helicoid $h$ from Lemma~\ref{le:unduloids}.
  The half-bubble $\plustilde B$ is the image under $h$
  of the rectangle $0<x<\pi$, $|y|<\tfrac14$, and we want
  to consider its Hopf image.

  The two $x$-coordinate lines in $\partial\plustilde B$
  with $y=\pm\tfrac14$ are \kHopf/ and project, as before, to the points
  $$
    p_\pm := \Pik\big(h(x,\pm\tfrac14)\big)
           = \k\,\cos\tfrac{n}2 \mp \j\,\sin\tfrac{n}2.
  $$
  The Hopf projection of~$\d\plustilde B$ is the minimizing great-circle
  arc~$e$ of length~$n$ (with~$\k=\Pik(h(0,0))$ as midpoint) connecting
  these points.  Indeed the $y$-coordinate lines with $x=0$ and $x=\pi$,
  which correspond to the necks of the unduloid,
  form (antipodal) arcs of length $n/2$ along the axis of the helicoid.
  Since this axis is $\i$-Hopf and $\i\perp\k$,
  these arcs project to the great-circle arc~$e$ with doubled length~$n$.
  The $y$-coordinate line with $x=\tfrac\pi2$, on the other hand,
  corresponds to the bulge of the unduloid, and lies along the polar
  axis of the helicoid.  This Hopf-projects to the complementary arc~$e'$
  of the same great circle (with~$-\k=\Pik(h(\tfrac\pi2,0))$ as midpoint).

  The $x$-coordinate lines foliate~$\plustilde B$ and
  are great semicircles along the rulings of the helicoid,
  corresponding to meridian arcs from neck to neck on the unduloid.
  Since the Hopf field of each is $u(y)\not=\k$,
  its \kHopf/ projection is a round circle~$c_y$ in~$\S^2$.
  The points $h(0,y)$, $h(\tfrac\pi2,y)$ and $h(\pi,y)=-h(0,y)$ are equally
  spaced along the ruling semicircle.
  The points $\pm h(0,y)$ project to the same point $q_y=e\cap c_y$,
  while $h(\tfrac\pi2,y)$ projects to the point $q'_y=e'\cap c_y$.
  Because of the equal spacing, $q_y$ and~$q'_y$ must be
  diametrically opposite in~$c_y$.
  As~$y$ increases from $-\tfrac14$ to $\tfrac14$,
  the points $q_y$ and~$q'_y$ each move at constant speed from $p_-$ to~$p_+$,
  along $e$ and~$e'$, respectively.
  It follows that the circles~$c_y$ foliate $\S^2\setm\{p_\pm\}$,
  and thus that $\Pik$ embeds $\plustilde B$ onto the slit sphere $\S^2\setm e$.
\end{proof}

\subsection{Spherical metrics on half-triunduloids}\label{ss:sphtri}

We saw how Hopf projection induces a spherical metric on a half-unduloid,
with each bubble corresponding to a slit sphere.
There is a similar decomposition of the spherical metric induced
on a half-triunduloid: it consists of slit spheres and one triangle,
corresponding to bubbles and the central piece of the triunduloid, respectively.
This decomposition will let us show that the spherical metric
depends only on the classifying triple of the triunduloid.

Lemma~\ref{le:slitray} characterizes the spherical metrics
on half-unduloids as lines of slit spheres; their
completion boundary develops to two distinct points in~$\S^2$.
Similarly, for triunduloids we consider
\emph{\threept/ metrics} which are spherical metrics on the disk whose
completion boundary consists of three points developing to a distinct triple.
The three points inherit a labeling from the triple; we require
that this be compatible with the orientation on the disk, meaning
that the labels increase counterclockwise.
\begin{lemma}\label{le:Phimetric}
  The metric induced on $M^+$ by the \kHopf/ projection of $\plustilde M$
  is a \threept/ metric, whose completion boundary develops to $\Psi(M)$.
\end{lemma}
\begin{proof}
  Let $D$ denote the disk $M^+$ with the spherical metric
  induced from Hopf projection of its cousin.
  Consider the horizontal boundary $M_0:= M\cap P$ of $M^+$,
  which has three connected components.
  The simply connected surface $M^+\cup M_0$ has a cousin with three boundary
  components consisting of Hopf fibers.
  Consider its quotient space~$D'$ obtained by identifying
  each boundary component to one point of some three-point set $T$.
  Observe that the Hopf projection descends to,
  and thus defines a spherical metric on, $D' = D \cup T$.  Moreover,
  $T$ Hopf projects to the triple~$\Psi(M)$.
  We claim $T$ is the completion boundary of $D$; that is, that $\hat D=D'$.

  We first show that $D'$ is complete,
  by decomposing it into a (nondisjoint) union of four complete sets.
  Choose open sets $E_i$ representing the ends of~$M$;
  removing them from $M^+\cup M_0$ leaves us with a compact set~$K$.
  Its quotient is a compact, and hence complete, subset of $D'$.
  On the other hand, according to Theorem~\ref{thkks}
  each end~$E_i$ is exponentially close to some unduloid~$U$ of necksize~$n_i$.
  Thus $\plustilde E_i$ is also exponentially close to~$\plustilde U$.
  Moreover, as shown in the proof of Theorem~\ref{th:necksizes},
  these two surfaces are bounded by rays on the same \kHopf/ circles.

  Choosing a smaller representative~$E_i$ if necessary,
  it follows that, when each is endowed with
  the spherical pullback metric from the cousin Hopf projection,
  $E_i^+$ is isometric to a subset of~$U^+$.
  Similarly, the closure of $E_i^+$, taken within $D'$,
  becomes isometric to a closed subset of the quotient of $U^+\cup U_0$.
  By Lemma~\ref{le:hopfund}, the latter quotient is
  the metric completion of a ray of spheres (with slit length~$n_i$);
  in particular it is complete.
  Therefore, the closures of the three ends~$E_i^+$ are
  complete subsets of $D'$.

  It remains to show that in the completion~$\hat D$
  the three points of~$T$ are indeed distinct.
  But the endpoints of the completion boundary arc
  of a ray of slit spheres have positive distance,
  and so each pair of points in $T$ also has positive distance.
\end{proof}

We say a spherical metric on the disk is embeddable if its developing map
is an embedding.  In particular, an \emph{embeddable triangle} is an
embeddable metric whose completion boundary is a loop developing onto
three minimizing great-circle arcs, connecting the points
of a triple $T\in\T$ pairwise.

\begin{lemma}\label{leimtr}
  Let the disk $D$ carry a \threept/ metric such that the
  completion boundary develops to a triple $T$,
  and let $\Gamma$ be the union of three minimizing
  great-circle arcs connecting the points in $T$ pairwise.
  Then $D$ is isometric to an embeddable triangle $\Delta$
  with a ray of slit spheres joined across each edge.
  The triangle $\Delta$ is isometric to the region in $\S^2$
  to the left of $\Gamma$.
\end{lemma}
\noindent
Note that in the case when $T$ includes a pair of antipodal points,
the loop $\Gamma$, and hence the triangle $\Delta$, is not uniquely
determined; for triunduloids, this happens exactly when there is
a cylindrical end.
In this case, we will agree to resolve the ambiguity by choosing $\Gamma$
to be a complete great circle, so that $\Delta$ has area~$2\pi$.
\begin{proof}
  Let $\phi\colon\hat D\to\S^2$ be the developing map.
  Any path in $\S^2\setm T$ lifts, by completeness, to some path in $D$.
  Thus if we remove the preimage $\phi^{-1}(T)$ from~$\hat D$,
  then $\phi$ restricts to a covering map.
\figtwo{The decomposition of the disk.}{%
  On the left, the open disk $D$ represents $M^+$.
  The Hopf projection of the cousin develops this disk onto~$\S^2$.
  The edges and vertices shown are the preimages of a minimizing
  triangle connecting the three vertices $\Psi(M)$ in~$\S^2$
  (shown in \figr{triunds}, bottom); they form the graph $G$
  which decomposes~$D$ into triangles.
  By ignoring the edges incident to the finite vertices,
  we can visualize the decomposition of Lemma~\ref{leimtr}:
  the disk~$D$ decomposes into a central triangle~$\Delta$,
  and three rays of slit spheres.
  On the right, we see the planar dual graph of $G$,
  covering a three-edge ``$\Theta$'' graph.
  The developing map of $D$ on the complement of the finite vertices
  is a covering map, homotopy equivalent to this graph covering.}

  Let us first deal with the case that the three arcs of $\Gamma$
  are disjoint from each other.
  The complement $\S^2\setm\Gamma$ then consists of two embeddable triangles.
  The completion $G$ of the preimage $\phi^{-1}(\Gamma)$ is a graph in $\hat D$.
  Its vertices are the three completion points and
  the points of $\phi^{-1}(T)\cap D$, which we call the finite vertices.
  The graph~$G$ has valence two at any finite vertex,
  since $\phi$ is an immersion there.

  The complement of $G$ consists of embeddable
  triangles isometric to those in $\S^2\setm\Gamma$.
  We claim there is a unique triangle $\Delta$ with no finite vertices.

  First, note that if some triangle had two (respectively, three)
  finite vertices, then it would be adjacent to the same triangle
  across all three of its edges.  The join of these two triangles
  would be a once-punctured sphere (respectively, an entire sphere).
  But neither of these is a three-point metric, and in neither case is
  there a way to join further triangles.

  Second, if every triangle had a single finite vertex,
  then these vertices would all develop to the same point of $T$,
  since this is true for adjacent triangles.
  But then our metric would be a line of slit spheres.  So
  there must be a triangle $\Delta$ with no finite vertices.  Its
  complement in $D$ consists of three rays of slit spheres, by
  Lemma~\ref{le:slitray}; in particular $\Delta$ is unique.  This
  proves the claim and establishes the desired decomposition.

  Finally, we must deal with the case where one arc $\overline{pr}$
  of~$\Gamma$ includes the third point $q$ of the triple in its
  interior.  In this case, we replace $\overline{pr}$ by its
  complementary arc in the same great circle.  The combinatorial
  argument above applies to this new $\Gamma'$, with the embeddable
  triangles replaced by domains isometric to hemispheres, whose
  completion boundaries develop to $\Gamma'$.  So there is a unique
  hemisphere whose three vertices form the completion boundary of $D$.
  We join this chosen hemisphere, across its nonminimizing edge, to
  the adjacent hemisphere.  This gives an embeddable triangle $\Delta$
  of area $4\pi$, whose completion boundary develops to the original
  $\Gamma$.
\end{proof}

For any triunduloid $M$, the last two lemmas find
within the spherical metric induced on~$M^+$
an embeddable triangle $\Delta$, which may be viewed as the subset of~$M^+$
enclosed by the three innermost necks.
The following proof that this metric is determined by~$\Psi(M)\in\T$
relies on the fact that this triangle is determined by $\Psi(M)$ alone.
\begin{proposition}\label{pr:trco}
  Let $M$ and $M'$ be two triunduloids having the same classifying triple
  $\Psi(M)=\Psi(M')\in\T$.  Then there is an isometry
  between the \threept/ metrics on $M^+$ and ${M'}^+$
  which preserves the labeling.
\end{proposition}
\begin{proof}
  We use Lemmas~\ref{le:Phimetric} and~\ref{leimtr} to decompose
  the spherical metric on the disk $\plustilde M$
  into an embeddable triangle $\Delta$ with three
  rays of spheres, and similarly, for $\plustilde{M'}$.

  The three vertices~$\Psi(M)$ uniquely determine the
  three minimizing arcs bounding the triangle~$\Delta\subset M^+$,
  except in the case when $\Psi(M)$ contains an antipodal pair;
  then we have agreed to choose $\Delta$ to have area $2\pi$.
  In any case, $\Delta$ is isometric to the triangle in~$\S^2$ which
  sits to the left of $\Pik(\d\wtilde\Delta)$ when traveling
  in the order of increasing label.

  Consequently, given $M$, $M'$ with $\Psi(M)=\Psi(M')\in\T$,
  the two embedded triangles $\Delta$ and $\Delta'$ are isometric.
  This isometry extends to the spherical metrics on the entire disks
  $M^+$ and ${M'}^+$
  by mapping corresponding rays of slit spheres onto one another.
\end{proof}

\subsection{A uniqueness theorem for triunduloid cousins}\label{ss:unique}

We are now ready to apply the maximum principle, in an appropriate bundle.
Given a triunduloid $M$, we use the \kHopf/ projection on the
cousin~$\plustilde M\subset\S^3$ to assign a spherical metric
to the disk $D=\plustilde M$ or $M^+$; this is the base manifold.
Now, $\S^3$ is the total space of the Hopf bundle over $\S^2$.
We can pull this bundle back
under the developing map to get a circle-bundle over $D$.
The cousin disk~$\plustilde M$ is transverse to the \kHopf/ fibers,
and can therefore be regarded as a section of this bundle.
To achieve a one-sided contact of two such sections, we
furthermore take the universal cover in the fiber direction,
obtaining a line bundle $E$ over the disk~$D$.  This covering exists,
because any circle bundle over a disk is a trivial bundle.

Note that at both steps we can pull back the Riemannian metric,
so that $E\to D$ becomes locally isometric to $\S^3\to\S^2$.
By Proposition~\ref{pr:trco} the metric on~$D$ depends only on the
classifying triple $T:=\Psi(M)$.
Hence the same is true of the pullback metric on $E=E(T)$.
We can then form the completion, $\hat E\to \hat D$; this adds
three $\k$-lines, over the three points in the completion boundary of $D$.
There is a fiberwise and isometric action of $\R$ on $\hat E$,
namely the pullback of the \kHopf/ flow on $\S^3$; we refer to this
as \emph{vertical translation} in $\hat E$.

For any triunduloid $M\in\M$ with $\Psi(M)=T$,
the cousin disk~$\plustilde{M}$ lifts to a minimal section of~$E(T)$.
Indeed, by Proposition~\ref{pr:phiimmerses}, $\plustilde{M}$ is transverse
to the \kHopf/ field and so locally lifts to a section of~$E(T)$;
but the bundle is trivial and so the lift is global.
Because cousins are defined only up to left translation in $\S^3$,
and we ignore rotations of $\S^2$,
this section is defined up to vertical translation.

Conversely, a minimal section~$\Sconj$ of~$E(T)$ is locally congruent
to a minimal surface in~$\S^3$, and so it has a \CMC/ cousin $S\subset\R^3$.
We are interested in the case when
the completion boundary of $\Sconj\subset\hat E$ contains
the three $\k$-lines over the completion boundary of~$D$.
(This means that the completion is not a section of~$\hat E$.)
Then the cousin disk~$S\subset\R^3$ has
three corresponding boundary components, each contained in a horizontal plane.
For a half-triunduloid, these planes would coincide, but in general,
their heights differ by the \emph{periods} of~$S$.
Let $\gamma$ be a (unit-speed) curve of finite length in~$S$,
connecting two of these boundary components.
By \eqn{firstH} the conormal vector $\etaconj$
to the cousin arc $\gammaconj$ in $\Sconj$ satisfies
$\etaconj=-\gammaconj\gamma'$.
Hence the period $P(\gamma)$, or signed height difference at the
endpoints of $\gamma$, is
\begin{equation}\label{eq:periods}
  P(\gamma)=\int_\gamma \langle \gamma',\k\rangle\,ds
  =-\int_\gammaconj \langle\gammaconj^{-1}\etaconj,\k\rangle\,ds
  =-\int_\gammaconj \langle\etaconj,\gammaconj\k\rangle\,ds,
\end{equation}
where we used the fact that left multiplication with $\gammaconj$
is an isometry.
\begin{lemma}\label{le:period}
  Consider two minimal sections~$\Sconj_1$,~$\Sconj_2$ of~$E(T)$ whose
  completions include the three completion fibers of $E$ as above.
  Suppose~$\gammaconj\subset\Sconj_1\cap\Sconj_2$ is a piecewise smooth
  curve of finite length,
  connecting two of the completion fibers,
  and suppose $\gammaconj$ has
  a one-sided neighborhood in~$\Sconj_1$
  which does not intersect $\Sconj_2$.
  If $\gamma_j$ denotes the cousin arc of $\gammaconj$ in $S_j\subset\R^3$,
  then the periods are unequal: $P(\gamma_1)\not=P(\gamma_2)$.
\end{lemma}
\begin{proof}
  Except at isolated points, two nonidentical minimal surfaces intersect
  along a curve transversely, that is, with different conormals.
  From our assumption that $\Sconj_1\cap \Sconj_2$ is empty to
  one side of~$\gammaconj$, and the fact that the surfaces are both
  transverse to the \kHopf/ field, we conclude that
  $\langle \etaconj_1-\etaconj_2,\,\gammaconj\k\rangle$
  has the same sign on a dense subset of $\gamma$.
  The claim then follows from~\eqn{periods}.
\end{proof}

Now we can apply the maximum principle to prove our uniqueness result.
We use the last lemma to deal with the noncompactness of the cousin disks.

\begin{theorem}\label{th:uniq}
  The classifying map $\Psi\colon\M\to\T$ is injective.
\end{theorem}
\begin{proof}
  Consider $M, N\in\M$ with $\Psi(M)=\Psi(N)=T\in\T$;
  let $n_1$, $n_2$, $n_3$ denote the common necksizes of~$M$ and~$N$.
  We regard the cousin disks~$\plustilde{M}$ and~$\plustilde N$ as sections
  of~$E(T)$.  Because $\d M^+$ consists of three curvature lines of
  infinite length in the horizontal plane, the completion boundary
  of $\plustilde M$ must include the three $\k$-lines added in $\hat E$;
  the same is true for $\plustilde N$.  Therefore, the lemma will apply.

  The base disk of the bundle $E(T)\to D$ consists, by Lemma~\ref{leimtr},
  of a central embeddable triangle and three rays of spheres.
  The bundle over the $i^{\text{th}}$ ray~$R_i$ is foliated by
  spherical helicoids, namely all vertical translates of the
  cousin of a half-unduloid with necksize $n_i$.
  By the asymptotics result of~\cite{kks} and Lemma~\ref{le:unduloids},
  each cousin is asymptotic, over~$R_i$, to one of these helicoid leaves.
  Thus, there is a well defined asymptotic height difference between
  $\plustilde M$ and~$\plustilde N$ at each end.

  If these three height differences are equal, then
  a vertical translation can be used to make~$\plustilde{M}$ asymptotic
  to~$\plustilde{N}$ over all three ends.
  We then claim $\plustilde{M}=\plustilde{N}$.
  These sections have the same smooth $\k$-lines as boundary,
  and $E(T)$ is foliated by minimal surfaces:
  the vertical translates of $\plustilde{M}$.
  We can find a leaf of the foliation which has a
  one-sided contact with~$\plustilde{N}$;
  the interior and boundary maximum principles then prove the claim.

  If the height differences are unequal, then after a vertical translation,
  $\plustilde M$ lies asymptotically above $\plustilde N$ on one end,
  but below it on some other end.  Let us then consider the
  open set $D'\subset D$ consisting of those points over which
  $\plustilde M$ lies strictly above~$\plustilde N$.
  By assumption, $D'$ has a component~$D''$ which includes one end.
  Then $\d D''\subset D$ is the projection of a piecewise smooth
  intersection curve $\gammaconj\subset \plustilde M \cap \plustilde N$,
  connecting two of the completion $\k$-lines.  This curve has finite length.
  Moreover, $\plustilde{M}$ and~$\plustilde{N}$ do not intersect over~$D''$.
  Thus we can apply Lemma~\ref{le:period} to find
  that the periods of $M^+$ and~$N^+$ cannot simultaneously vanish.
  This contradicts the fact that both $M^+$ and $N^+$ are half triunduloids.
\end{proof}

\section{Properness}\label{se:prop}

Korevaar and Kusner~\cite{kk} prove that each embedded $k$-unduloid of
genus~$g$ is enclosed in a union of solid cylinders and spheres of fixed
radii, the number of which depends only on~$g$ and $k$.
Making use of this enclosure theorem they derive area and curvature estimates.
We want to get similar estimates for triunduloids,
which may only be \alex-emb/.

The area estimates and enclosure theorem from~\cite{kk}, and their proofs,
carry over literally to this setting.
\begin{theorem}[\cite{kk}]\label{th:kk}
There is a constant $C$ such that for any triunduloid $M$,
we have $\area(M\cap B_r)\le Cr^2$ for any euclidean ball $B_r$ of radius $r$.
Furthermore, there is a ball $B(M)$ of diameter $21$ such that
$M\setm B(M)$ has three noncompact connected components, namely the three
ends of $M$; each of these is contained in a
solid cylinder of radius $3$ about some ray.
\end{theorem}
The curvature estimate, on the other hand, requires a new proof.
One step of the argument in~\cite{kk}
relies on the theory of embedded minimal surfaces of
finite total curvature.  Here, we use instead certain special
properties of Alexandrov-symmetric minimal surfaces.

We write $|A|^2$ for the sum of the squares of the principal curvatures
of $M$, that is, $|A|$ is the pointwise norm of
the second fundamental form of $M$.

\begin{lemma}\label{le:minimallimit}
Suppose $\sup_{M_k}|A|\to\infty$ for a sequence $M_k$ of triunduloids.
Then there are points $p_k\in M_k$ such that
the rescaled surfaces $M'_k := |A(p_k)|\, (M_k-p_k)$ have a subsequence
converging to an Alexandrov-symmetric (nonplanar) genus-zero
minimal surface $N$ of finite total curvature.
\end{lemma}
\noindent
Here, convergence is in the sense of local graphs.
The limit can have positive multiplicity,
but this must be finite by the area bounds.
\begin{proof}
  We choose the points $p_k\in M_k^+$
  such that on $M'_k$ we have $|A|\le 2$ on balls $M'_k\cap B_{R_k}(0)$,
  whose radii satisfy $|A(p_k)|\ge R_k\to \infty$ (see~\cite[p.~487]{kks}).
  Due to this uniform curvature bound and the area bound, the $M_k'$
  subconverge to a minimal surface~$N$, together with their normals.
  On~$N$, we have $|A(0)|=1$, so $N$ is not a plane.

  To bound the total curvature of $N$, we follow~\cite{kks} in
  choosing connected subsets
  $S_k\subset M_k$ with uniformly bounded area such that
  $M_k\cap B_1(p_k)\subset S_k$ and $\d S_k$ consists of geodesic loops
  (around necks of the surface).
  Since the~$M_k$ are triunduloids, the number of loops needed is at most three;
  the Gauss-Bonnet theorem then uniformly bounds the total curvature of~$S_k$.
  Because total curvature is invariant under rescaling,
  and the rescaled~$S_k$ converge to~$N$,
  we find that $N$ has finite total curvature $\int_N |A|^2$.
  (See \cite[\paragraph 5]{kk} and~\cite[p.486--488]{kks} for details.)

  The~$M_k$ are Alexandrov-symmetric with respect to the horizontal plane $P$.
  Suppose the symmetry planes of the $M'_k$ are not at bounded distance
  from the origin.  Alexandrov symmetry implies that on ${M'_k}^+$
  the Gauss map takes values in the upper hemisphere.
  Thus by convergence of normals, the minimal surface $N$ also has Gauss image
  in this hemisphere, contradicting the fact that $N$ is nonplanar.
  So instead, the $M'_k$ have symmetry planes at bounded
  distance from the origin; a subsequence will have convergent symmetry planes.
  Convergence of normals guarantees that
  $N$~has Alexandrov symmetry in the limit plane.
\end{proof}

To make use of this lemma, we now consider
minimal surfaces of finite total curvature.  Each end has a
well defined normal and winding number~\cite[Sect.~2.1]{hk}.  If an end
has winding number one, it is embedded and it is asymptotically
\emph{catenoidal} or \emph{planar}.
We call an end \emph{horizontal} when
its limiting normal is vertical, and vice versa.

We also need to recall the notion of forces \cite[\paragraph 3]{kks}.
If~$\Gamma$ is a cycle on a surface~$M$ with constant mean curvature~$H$,
the vector $f(\Gamma):= \int_{\Gamma}\eta + 2H\int_{K}\nu \in\R^3$ is
called the \emph{force} of~$\Gamma$.  Here $\eta$ is the
conormal along~$\Gamma$ and $K$~is a 2-chain in $\R^3$ with boundary~$\Gamma$
and normal~$\nu$.  The first variation formula for area shows that the
force depends only on the homology class~$[\Gamma]\in H_1(M)$.

Note that each annular end $E$ of $M$ defines a cycle,
and thus has an associated force.
If $E$ is exponentially asymptotic to a surface of revolution
(like a catenoid or unduloid)
its force equals that of the asymptotic surface, and in particular
is directed along the axis with magnitude
\begin{equation}\label{eq:weight}
  n\big(1-H\tfrac{n}{2\pi}\big),
\end{equation}
where $n$ is the necksize, or length of
the shortest closed geodesic.
This magnitude is nonzero for catenoids ($H=0$, $n>0$) and for unduloids
($H=1$, $0<n\le\pi$).

\begin{lemma}\label{lemsce}
  Let $N$ be a minimal surface of finite total curvature, having
  Alexandrov symmetry with respect to a horizontal plane~$P$.
  Then each end is either horizontal or vertical, and each vertical
  end is catenoidal.
  Furthermore, if $N\cap P$ contains an unbounded curve, then there is a
  simple closed geodesic on~$N$, symmetric with respect to~$P$,
  which has nonzero force.
\end{lemma}

\begin{proof}
  An Alexandrov-symmetric surface meets the horizontal symmetry plane~$P$
  perpendicularly.  An end not meeting $P$ must be horizontal; an end
  meeting~$P$ is vertical and is itself Alexandrov-symmetric.
  At any planar end of a minimal surface,
  the Gauss map has a branch point~\cite{hk},
  so a vertical planar end contradicts Alexandrov symmetry.
  Similarly, a vertical end with higher winding number cannot be decomposed
  into symmetric halves on the two sides of $P$, so it cannot
  be Alexandrov-symmetric.  Thus a vertical end is catenoidal.

  If $N\cap P$ is not a union of compact curves,
  there must be an end~$E$ meeting~$P$.  By the above, $E$ is vertical
  and catenoidal, so it has nonzero force.
  Consider a simple closed loop around~$E$, which is nontrivial in~$N$.
  Because $N$ has negative Gauss curvature except at isolated points,
  there is a unique simple closed geodesic $\gamma$ homotopic to this loop.
  The reflection of~$\gamma$ in~$P$ is also a geodesic,
  and so coincides with $\gamma$.  Thus $\gamma$ is symmetric;
  its force equals that of~$E$.
\end{proof}

\begin{proposition}\label{pr:curvbound}
  For each $\epsilon>0$, there exists a constant $C$ such that
  any triunduloid~$M$ with necksizes $n_1,n_2,n_3>\epsilon$
  has uniformly bounded curvature, $|A|\le C$.
\end{proposition}
\begin{proof}
  Suppose we have $\sup_{M_k}|A|\to\infty$
  for some sequence of triunduloids~$M_k$
  with necksizes greater than~$\epsilon$.
  Then the rescaled surfaces $M'_k$, as in Lemma~\ref{le:minimallimit},
  converge to an Alexandrov-symmetric minimal surface~$N$.
  The~$M_k$ have genus~$0$, so there are no closed loops in $M_k\cap P$.
  Since the $M'_k$ converge to $N$ with finite multiplicity,
  there are no closed loops in~$N\cap P$ either.
  But $N$ meets $P$, and hence does so in an unbounded curve.

  Thus, we can apply Lemma~\ref{lemsce} to get a
  closed geodesic $\gamma$ in $N$ with nonzero force.
  This is the limit of closed curves $\gamma'_k$ on $M'_k$,
  whose lengths and forces converge to those of $\gamma$.
  So on the unrescaled $M_k$, the corresponding curves $\gamma_k$
  have nonzero forces, which converge to zero.
  Because $\gamma_k$ has nonzero force, it is nontrivial; a simple nontrivial
  curve on a triunduloid is homologous to the boundary of one of the ends.
  By the necksize bound, the force of \eqn{weight} is bounded away from zero,
  a contradiction.
\end{proof}

\begin{theorem}\label{th:proper}
  The classifying map $\Psi\colon\M\to\T$ is continuous and proper.
\end{theorem}
\begin{proof}
  Consider a compact subset $\mathcal K\subset\T$.
  By Theorem~\ref{th:necksizes} there is an $\epsilon>0$ such that each surface
  $M\in\Psi^{-1}(\mathcal K)$ has necksizes bounded below by~$\epsilon$,
  and thus uniformly bounded curvature by Proposition~\ref{pr:curvbound}.

  Suppose two triunduloids in $\Psi^{-1}(\mathcal K)$ are sufficiently close
  in Hausdorff distance on an open domain in $\R^3$ with compact closure.
  Because of the uniform curvature bounds,
  either one can then be written as a normal graph over the other
  in any closed subset of the domain.  By elliptic regularity,
  these subsets of the triunduloids are then also $C^1$ (or~$C^\infty$) close.
  Thus, up to left translation in~$\S^3$,
  the corresponding subsets of the cousins are again close,
  and so are the Hopf projections of their boundaries.
  We conclude that $\Psi$ is continuous.

  To show that $\Psi^{-1}(\mathcal K)$ is compact,
  let us now prove that any sequence $M_k\in \Psi^{-1}(\mathcal K)$
  subconverges with respect to the topology on~$\M$.
  The~$M_k$ are given only up to horizontal translation; to
  get convergence, we choose representatives such that
  the enclosing balls $B(M_k)$ from Theorem~\ref{th:kk}
  are centered at the origin.
  Due to the area and curvature estimates we find a subsequence of the~$M_k$
  that converges on each compact subset to some \CMC/ surface~$M$.

  By smooth convergence, it is clear that $M$ is \alex-emb/.
  Because a sequence of contractible loops has a contractible limit,
  $M$ has genus zero and at most three ends.
  On the other hand, by the enclosure theorem, on each $M_k$
  there are three closed curves in $B_{12}(0) \setm B_{11}(0)$
  which bound the three ends.  These three sequences of
  curves subconverge to three closed curves which are disjoint in
  the limit surface $M$, and thus bound three different ends.
  Thus $M$ is a triunduloid, as desired.
\end{proof}

\section{Surjectivity}\label{se:surj}

In this section we show that $\Psi$ maps $\M$ onto the three-ball~$\T$.
We depend on the fact that the moduli space~$\M$ of triunduloids 
is locally a real analytic variety.
We show that this variety has dimension three by
establishing the existence of a nondegenerate triunduloid.
Between manifolds of the same dimension, a proper injective map must 
also be surjective, and in fact a homeomorphism.  
The structure theory for real-analytic varieties lets us extend this
result to the classifying map $\Psi\colon\M\to\T$.

\subsection{The real-analytic variety structure of $\M$}

The moduli space of \CMC/ surfaces of any fixed finite topology
is locally a real-analytic variety, by the local structure theory
of Kusner, Mazzeo, and Pollack~\cite{kmp}.  In order to discuss the regular
points of this variety, we recall the notion of nondegeneracy:

\begin{definition}\label{dedege}
  A minimal or \textsc{cmc} surface $M$ is called \emph{degenerate}
  if the Jacobi equation $\Delta_{{}_M} u+|A|^2 u=0$
  has a nonzero square-integrable solution $u\in L^2(M)$.
  Otherwise $M$ is called \emph{nondegenerate}.
\end{definition}

An implicit function theorem argument~\cite{kmp} shows that, in general,
the moduli space of \CMC/ surfaces of fixed genus with $k$~embedded
ends is a finite-dimensional real-analytic variety, and that
near a nondegenerate \CMC/ surface, this variety is a $(3k-6)$--manifold.
In particular, for triunduloids:
\begin{theorem}[\cite{kmp}]\label{th:kmp}
  The moduli space $\M$ is locally
  a real-analytic variety of finite dimension.
  In a neighborhood of any nondegenerate triunduloid,
  $\M$ is a manifold of dimension three.
\end{theorem}

We remark that the converse is not necessarily true: 
a three-manifold point of~$\M$ may still be a degenerate triunduloid.
Therefore, while our Main Theorem says that $\M$ is a three-manifold
everywhere, it remains an open question whether all triunduloids are
nondegenerate.  

A \CMC/ moduli space can have dimension different from
the expected~$3k-6$, and thus consist entirely of degenerate surfaces.
In fact, the moduli space of \CMC/ surfaces with a fixed finite topology
can have components of arbitrarily large dimension; examples include the
so-called multibubbletons with two ends, and the \CMC/ tori.

On the other hand, it seems likely that most \alex-emb/ surfaces,
and perhaps all the Alexandrov-symmetric ones, are nondegenerate.
While the methods of~\cite{kgb} give certain degenerate
embedded triply periodic \CMC/ surfaces (considered as
embeddings of a compact surface into a flat three-torus), no similar example
is known with finite topology.

\subsection{Nondegenerate triunduloids}

The existence of the injective map $\Psi\takes \M\to\T$
shows that the variety $\M$ cannot have dimension greater than three.
However, to apply topological methods to get surjectivity, we need
to know the dimension is exactly three.  To get this, we want to
prove the existence of a nondegenerate triunduloid, near which $\M$
will be a three-manifold.

Montiel and Ros~\cite{mr} have a theorem about nondegeneracy of minimal
surfaces.  Here we specialize this result to
\emph{minimal $k$-noids}, minimal surfaces of finite total curvature and
genus zero with only catenoidal ends.
\begin{theorem}[{\cite[Cor.~15]{mr}}]
  Suppose $M$ is a minimal $k$-noid
  such that the Gauss image of its umbilic points
  is contained in a great circle in~$\S^2$.
  Then the only bounded Jacobi fields on~$M$ are those induced by
  translations in~$\R^3$.
\end{theorem}
\noindent
In our reformulation we have used the fact that the Gauss map of
a $k$-noid does not branch at the ends.
\begin{corollary}\label{co:trinoids}
  All minimal trinoids are nondegenerate.
\end{corollary}
\begin{proof}
Note that a trinoid $M$ has total curvature $-8\pi$.
Hence the Gauss map has degree two and so by the Riemann-Hurwitz formula
at most two branch points, which are the umbilic points of $M$.
Therefore the theorem applies.  If $M$ were degenerate, there
would be an $L^2$~Jacobi field, which would be bounded, and thus
induced by a translation.  But a translation of a catenoidal end
(even if orthogonal to the axis) has normal component which fails to be~$L^2$.
\end{proof}

Mazzeo and Pacard~\cite{mp} find \CMC/ surfaces by gluing
unduloids or nodoids of small necksize to
the ends of an appropriately scaled nondegenerate minimal $k$-noid.
\begin{theorem}[\cite{mp}]\label{th:mp}
  Let $M_0\subset\R^3$ be a nondegenerate minimal $k$-noid that
  can be oriented so that each end has normal pointing inwards towards its
  asymptotic axis.
  Then for some $\epsilon_0$,
  there is a family $M_{\epsilon}$, $0<\epsilon<\epsilon_0$,
  of nondegenerate \CMC/ surfaces with
  $k$~embedded ends, such that on any compact set in~$\R^3$,
  the dilated surfaces
  $\tfrac 1{\epsilon} M_{\epsilon}$ converge smoothly and uniformly to~$M_0$.
\end{theorem}
Gluing an unduloid to each end of a minimal trinoid
yields a surface which, though not embedded for small necksizes,
is \alex-emb/.
Thus we can combine Corollary~\ref{co:trinoids}
with Theorem~\ref{th:mp} to yield:
\begin{corollary}\label{co:nondegtriund}
  There exists a nondegenerate triunduloid.
\end{corollary}

\subsection{Topological arguments for surjectivity}

We recall (see~\cite[\paragraph 7-4]{mun}) that a space $Y$ is
\emph{compactly generated} provided that any subset $A$ is open (resp.~closed)
in $Y$ if and only if $A\cap K$ is open (resp.~closed) in $K$
for every compact $K\subset Y$.
Any locally compact space $Y$ is compactly generated, as is any
space $Y$ which has a countable base at every point.

\begin{lemma}\label{le:homeo}
If $f\takes X\to Y$ is a continuous, proper, injective map from
an arbitrary space $X$ to a compactly generated Hausdorff space $Y$,
then $f$ is an embedding, that is,
a homeomorphism onto its image $f(X)$.
Furthermore, this image is closed in $Y$.
\end{lemma}
\begin{proof}
Because $f$ is a continuous injection, it is an embedding
if and only if it is a closed (or open) map onto its image.
Given a closed set $A\subset X$, we claim $f(A)$ is closed in~$Y$,
and hence in $f(X)$; applying the claim to $A=X$ will confirm the
last statement of the lemma.

We check the claim
by showing $f(A)\cap K$ is closed for all compact $K\subset Y$.
Since $f$ is proper, $f^{-1}(K)$ is compact, and so the closed subset
$A\cap f^{-1}(K)$ is also compact.  Thus its image $f(A)\cap K$
is compact, and (since $Y$ is Hausdorff) closed.
\end{proof}

\begin{definition}
  Let $X$ be a locally finite $d$-dimensional simplicial complex
  and $S$ a $(d{-}1)$-simplex in $X$.  The \emph{valence} of~$S$ is the
  number of $d$-simplices which contain~$S$ as a face.
  We say the complex~$X$ is \emph{borderless}
  if no $(d{-}1)$-simplex has valence $1$.
\end{definition}

\begin{theorem}\label{th:homeo}
  If $f\takes X\to Y$ is a continuous, proper, injective map from
  a borderless $d$-complex $X$ to a connected $d$-manifold $Y$,
  then $f$ is surjective and thus (by the previous lemma) a homeomorphism.
\end{theorem}
\begin{proof}
  Consider the union of the closed $d$-simplices in~$X$ and remove its
  $(d{-}2)$-skeleton.   Let~$X'$ be the closure (in~$X$) of some connected
  component of this set.  Its $(d{-}1)$-simplices have
  valence at least two because $X$ is borderless.  Since $X'$~is closed,
  the restriction $f|_{X'}$ is still proper.
  We will show this restriction is surjective.
  (This implies that $X=X'$ and that every valence is exactly two.)

  By Lemma~\ref{le:homeo}, $f$ is an embedding and $f(X')$ is closed in~$Y$.
  In particular, the image of the $(d{-}2)$-skeleton $X'_{d-2}$ is
  a properly embedded $(d{-}2)$-complex in the connected $d$-manifold $Y$;
  therefore its complement $B$ is connected.
  Letting $Z$ be the complement of $X'_{d-2}$ in $X'$,
  we have that $f(Z)=B\cap f(X')$ is a closed subset of $B$.
  We claim that $f(Z)$ is also open in~$B$.
  It then follows by the connectedness of~$B$ that $f$ maps $Z$ onto~$B$.
  Since $\bar Z = X'$ and $f(X')$ is closed,
  $f(X') \supset \bar B = Y$, as desired.

  To prove the claim, suppose first that $z\in Z$ is interior to a
  $d$-simplex.  Since the dimension of~$Y$ is also~$d$, and $f$
  is an embedding, $f(Z)$ contains an open neighborhood of~$f(z)$.
  Otherwise $z$ is interior to a $(d{-}1)$-simplex $S$.
  Because $X'$ is borderless, $z$~is in the interior of the union of two
  $d$-simplices with common face~$S$,
  so again $f(Z)$ contains a neighborhood of~$f(z)$.
\end{proof}

\subsection{Structure of real-analytic varieties and proof of the Main Theorem}

To apply our topological theorem to real-analytic varieties, we recall the 
notion of an Euler space (see~\cite[Def.~4.1]{fm}).

\begin{definition}
  Let $X$ be a $d$-dimensional locally finite simplicial complex.
  We say $X$ is a \emph{$\mod2$--Euler space} if the link of each
  simplex in $X$ has even Euler characteristic, or, equivalently, if
  the local Euler characteristic $\chi(X,X\setm x)$ is odd for all $x\in X$.
\end{definition}
Recall that the link of a $(d{-}1)$-simplex $S \subset X$ is
obtained as follows: consider the vertices of the $d$-simplices with
$S$ as a face; those which are not in~$S$ itself form the link of~$S$.
The valence of~$S$ gives the number of points in the link.  For a
$\mod2$--Euler space the link of any $(d{-}1)$-simplex consists of an
even number of vertices.  Therefore, a $\mod2$--Euler space is
borderless.

{\L}ojasiewicz~\cite{loj} showed real-analytic varieties are triangulable
as locally finite simplicial complexes,
and Sullivan discovered that they have the $\mod2$--Euler property
(see \cite[Thm.~4.4]{fm} or~\cite{sul,bv,har1}).
As a consequence of these results, we have:
\begin{proposition}\label{pr:borderless}
  A $d$-dimensional real analytic variety, or any space locally
  homeomorphic to one, is a borderless $d$-complex.
\end{proposition}

Making use of our previous results on injectivity and properness we can
now prove:
\begin{theorem}
  The classifying map $\Psi\colon\M\to \T$ is a homeomorphism.
\end{theorem}
\begin{proof}
  To apply Theorem~\ref{th:homeo} to $\Psi\colon\M\to \T$, we need to
  verify its assumptions.

  The space $\T$ of triples is clearly a connected three-manifold.  The
  classifying map $\Psi$ is proper and injective by
  Theorems~\ref{th:proper} and~\ref{th:uniq}.  From the first part of
  Theorem~\ref{th:kmp}, the moduli space~$\M$ of triunduloids is
  locally a real analytic variety; injectivity of $\Psi$ shows it has
  dimension at most three.  So Proposition~\ref{pr:borderless}
  implies~$\M$ is a borderless $d$-complex for some $d\le3$.

  It remains to check that $\M$ has dimension~$d=3$.  By
  Corollary~\ref{co:nondegtriund} there is a nondegenerate
  triunduloid, so by the second part of Theorem~\ref{th:kmp}, $\M$ is
  a three-manifold in some neighborhood of this triunduloid.
  Theorem~\ref{th:homeo} now applies to conclude $\Psi$ is a homeomorphism.
\end{proof}
Together with Theorem~\ref{th:necksizes} this completes the
proof of our Main Theorem.


\end{document}